\documentclass{amsart}

\usepackage{newtxtext}
\usepackage{newtxmath}
\usepackage{graphicx}
\usepackage{multicol}
\usepackage{footmisc}
\usepackage[normalem]{ulem}
\usepackage[x11names]{xcolor}

\usepackage{comment}

\usepackage{graphics}
\usepackage{url}
\usepackage{amscd}
\usepackage{textcomp}
\usepackage{enumitem}
\usepackage[all]{xy}
\usepackage{mathtools} 

\usepackage{tikz}
\usepackage{ifthen}
\usetikzlibrary{shapes.misc, arrows}
\tikzset{cross/.style={cross out, draw=black, minimum size=2*(#1-\pgflinewidth), inner sep=0pt, outer sep=0pt},
cross/.default={1pt}}

\definecolor{darkblue}{rgb}{0,0,0.4} 
\usepackage{xr-hyper}
\usepackage[colorlinks=true, citecolor=darkblue, filecolor=darkblue, linkcolor=darkblue,urlcolor=darkblue]{hyperref} 
\usepackage[all]{hypcap}

\usepackage{tikz}
\usetikzlibrary{arrows.meta,decorations.markings,arrows}
\usepackage{caption}

\newtheorem{theorem}{Theorem}[section]

\newtheorem{lemma}[theorem]{Lemma}

\theoremstyle{definition}

\newtheorem{definition}[theorem]{Definition}
\newtheorem{exercise}[theorem]{Exercise}
\theoremstyle{remark}
\newtheorem{remark}[theorem]{Remark}
\newtheorem{exampl3}[theorem]{Example}
\numberwithin{equation}{section}

\newcommand{\grid}{{\mathbb {G}}}
\newcommand{\domain}{\mathcal {D}}

\newcommand{\mc}{\mathcal}

\newcommand{\eps}{\varepsilon}

\newcommand{\set}[1]{\left\{#1\right\}}

\newcommand{\del}{\partial}

\newcommand{\Z}{\mathbb{Z}}

\newcommand{\R}{\mathbb{R}}
\newcommand{\C}{\mathbb{C}}
\newcommand{\CP}{\mathbb{CP}}
\newcommand{\T}{\mathbb{T}}

\newcommand{\Cat}{\mathcal{C}}   
\DeclareMathOperator{\Ob}{Ob}    
\DeclareMathOperator{\TopSpaces}{Top}  
\DeclareMathOperator{\Ab}{Ab}    
\newcommand{\Mod}{\mathcal{M}}    
\newcommand{\Sym}{\mathrm{Sym}} 

\DeclareMathOperator{\Crit}{Crit}
\newcommand{\Ccell}{\widetilde{\mathcal{C}}_*^{\textrm{cell}}}
\newcommand{\CMorse}{\widetilde{\mathcal{C}}_*^{\textrm{Morse}}}

\newcommand{\x}{\mathbf{x}}
\newcommand{\y}{\mathbf{y}}
\newcommand{\z}{\mathbf{z}}
\newcommand{\w}{\mathbf{w}}

\DeclareMathOperator{\CKh}{CKh}

\definecolor{sgreen}{cmyk}{1,0,1,0.6}

\newcommand{\drawO}[2]{\draw[thick] (#1-0.5,#2-0.5) circle (0.3);}
\newcommand{\drawX}[2]{
    \path (#1-0.15,#2-0.15) node(localx) {} 
        (#1-0.85,#2-0.85) node(localy) {};
    \draw[thick] (localx) -- (localy);
    \path (#1-0.15, #2-0.85) node(localx) {} 
        (#1-0.85,#2-0.15) node(localy) {};
    \draw[thick] (localx) -- (localy);
    }
    
\newcommand{\SARdoublybroken}[3]{
\draw (#1#2#3) + (0,0.75) node (#1#2#3-TT){};
\draw (#1#2#3) + (0,0.25) node (#1#2#3-T){};
\draw (#1#2#3) + (0,-0.25) node (#1#2#3-B){};
\draw (#1#2#3) + (0,-0.75) node (#1#2#3-BB){};
\draw (#1#2#3-TT.center) edge[out=-150, in=150] (#1#2#3-T.center);
\draw (#1#2#3-TT.center) edge[out=-30, in=30] (#1#2#3-T.center);
\draw (#1#2#3.center) + (0,0.5) node {\tiny $#1$};
\draw (#1#2#3-T.center) edge[out=-150, in=150] (#1#2#3-B.center);
\draw (#1#2#3-T.center) edge[out=-30, in=30] (#1#2#3-B.center);
\draw (#1#2#3.center) node {\tiny $#2$};
\draw (#1#2#3-B.center) edge[out=-150, in=150] (#1#2#3-BB.center);
\draw (#1#2#3-B.center) edge[out=-30, in=30] (#1#2#3-BB.center);
\draw (#1#2#3.center) + (0,-0.5) node {\tiny $#3$};
}

\newcommand{\SARbrokendown}[3]{
\draw (#3d) + (0,0.75) node (#3d-T){};
\draw (#3d) + (0,-0.5) node (#3d-B){};
\draw (#3d-T.center) edge[out=-160, in=160, looseness=1.7] (#3d.center);
\draw (#3d-T.center) edge[out=-20, in=20, looseness=1.7] (#3d.center);
\draw (#3d.center) + (0,0.375) node {\tiny $#1+#2$};
\draw (#3d.center) edge[out=-150, in=150] (#3d-B.center);
\draw (#3d.center) edge[out=-30, in=30] (#3d-B.center);
\draw (#3d.center) + (0,-0.25) node {\tiny $#3$};
}

\newcommand{\SARbrokenup}[3]{
\draw (#1u) + (0,0.5) node (#1u-T){};
\draw (#1u) + (0,-0.75) node (#1u-B){};
\draw (#1u-T.center) edge[out=-150, in=150] (#1u.center);
\draw (#1u-T.center) edge[out=-30, in=30] (#1u.center);
\draw (#1u.center) + (0,0.25) node {\tiny $#1$};
\draw (#1u.center) edge[out=-160, in=160, looseness=1.7] (#1u-B.center);
\draw (#1u.center) edge[out=-20, in=20, looseness=1.7] (#1u-B.center);
\draw (#1u.center) + (0,-0.375) node {\tiny $#2+#3$};
}

\begin{document}

\title{Spectra in Khovanov and knot Floer theories}

\author{Marco Marengon}%
\address{HUN-REN Alfr\'ed R\'enyi Institute of Mathematics, Budapest, Hungary}%
\email{\href{mailto:marengon@renyi.hu}{marengon@renyi.hu}}%

\author{Sucharit Sarkar}
\address{UCLA, Los Angeles, CA}%
\email{\href{mailto:sucharit@math.ucla.edu}{sucharit@math.ucla.edu}}%

\author{Andr\'as Stipsicz}%
\address{HUN-REN Alfr\'ed R\'enyi Institute of Mathematics, Budapest, Hungary}%
\email{\href{mailto:stipsicz@renyi.hu}{stipsicz@renyi.hu}}%

\date{}

\begin{abstract}
These notes provide an introduction to the stable homotopy types in Khovanov theory (due to Lipshitz-Sarkar) and in knot Floer theory (due to Manolescu-Sarkar).
They were written following a lecture series given by Sucharit Sarkar at the R\'enyi Institute during a special semester on ``Singularities and low-dimensional topology'', organised by the Erd\H os Center.
\end{abstract}


\maketitle

\tableofcontents


\section{Introduction}

Invariants of knots and links witnessed an unparalleled development in the past couple of decades. 
The subject of knot theory was born at the end of the XIX century, while a systematic mathematical study began
with work of Reidemeister in the 1920s \cite{Reidemeister}. Alexander's work
(and in particular the introduction of the Alexander
polynomial invariant \cite{Alexander}) in the late 1920s was a remarkable
step in understanding this subject. The next major 
development was brought by work of Vaughan Jones in the 
80s: his introduction of the Jones polynomial \cite{Jones}, and 
subsequent generalizations, the HOMFLY-PT polynomial and 
the definition of quantum invariants, reconnected the
subject with theoretical physics. 

After the turn of the
century a new paradigm, \emph{categorification}, provided
deep insight into the (still mysterious) world of knots.
Khovanov's definition of his homologies \cite{K:Khovanov} refined the 
Jones polynomial, while work of Ozsv\'ath and Szab\'o \cite{OSz:HFK},
and independently of Rasmussen \cite{R:HFK} produced a family of homology theories having the Alexander polynomial as their (graded)
Euler characteristic. It was a routine exercise to find 
knots which are not distiguished by their Jones (or Alexander)
polynomials, but have distinct Khovanov (or knot Floer)
homologies. Indeed, both Khovanov and Heegaard Floer
homologies are \emph{unknot detectors} (i.e. the 
unknot is characterized by its homology in either theory).
For comparison, this statement is false for the Alexander
polynomial, and an open problem for the Jones polynomial.

Recently a further level of abstraction, \emph{spacification}, brought new tools
in knot theory. In this context one wants 
to construct topological spaces (or, more generally, spectra) with the property that these objects are knot invariants, 
and their singular homology (in the appropriate sense) 
is the Khovanov homology (or knot Floer homology) of the
underlying knot. Having a spectrum provides the opportunity of applying further structures (such as Steenrod operations) and possibly other functors (like $K$-theory) to find further knot invariants. The program has been completed for Khovanov homology \cite{LS:KhovanovHomotopyType}, and the first step has been taken for knot Floer homology \cite{ManolescuSarkar}.

These notes are based on lectures of Sucharit Sarkar given at a
focused week of the special semester \emph{Singularities and
  low-dimensional topology} at the Erd\H{o}s Center, Budapest in the
Spring semester of 2023.  We thank Michael Willis for his comments on
an earlier version of these notes.

MM acknowledges that: This project has received funding from the
European Union's Horizon 2020 research and innovation programme under
the Marie Sk{\l}odowska-Curie grant agreement No.\ 893282. AS was
partially supported by NKFIH Grant K146401.

\section{Framed flow categories}

The basic idea of \emph{spacification} via flow categories is to start from a chain complex $C$ together with a given basis and construct a `space' (more precisely a CW spectrum) $X$ with the property that its reduced cellular chain complex $\Ccell(X)$ is isomorphic to $C$. The following remark explains why we really need to choose a basis of $C$ for our construction.

\begin{remark}
\label{rem:spac-basis}
If we have a chain complex $C$, which we assume to be finitely generated for simplicity, and we disregard its basis and grading, then we can always make a basis change so that
\[
C \cong \Z^k \oplus \bigoplus_i (\Z \xrightarrow{k_i} \Z),
\]
where the boundary map on $\Z^k$ is zero and the symbol $\Z \xrightarrow{k_i} \Z$ denotes a two-step complex where the boundary map is the multiplication by $k_i$,
and from here we can construct a CW complex $M(C)$ with $\Ccell(X)=C$ simply by taking an appropriate wedge of spheres (for the free part) and other simple spaces (for the 2-step complexes) as described in \cite[Example 2.40]{Hatcher}.
The space $M(C)$ is often called a \emph{Moore space}, and its construction can be refined to take gradings in account as well.

The drawback of Moore spaces is that they are not (and they cannot be made) functorial. More precisely, by a result of Carlsson \cite{C:Moore} the functors
\[
H_*, H_n \colon \TopSpaces_* \to \Ab
\]
do not admit functorial sections $M, M_n \colon \Ab \to \TopSpaces_*$.

This is in sharp contrast with the case of homotopy groups, where the functor
\[
\pi_n \colon \TopSpaces_* \to \Ab,
\]
defined for $n>1$, does admit a functorial section, given by the Eilenberg-Maclane spaces.
\end{remark}

\subsection{A recap of Morse theory}
Let $f \colon M \to \R$ be a proper Morse-Smale function with finitely many critical points on a Riemannian manifold $M$.
Each index-$k$ critical point gives rise to a \emph{downward disc} $\mathbb{D}^k$, which together with all other downward discs presents $M$ as a (based) CW complex.

\begin{figure}
\begin{center}
    \begin{tikzpicture}
        \def\x{1};
        \def\height{4};
        
        \draw [dashed, gray]  (-2*\x,3*\x) arc[start angle=180,end angle=0,x radius=0.5*\x,y radius=0.25*\x]; 
        \draw[gray] (-2*\x,3*\x) arc[start angle=180,end angle=360,x radius=0.5*\x,y radius=0.25*\x]; 
        \draw [dashed, gray]  (-2*\x,\x) arc[start angle=180,end angle=0,x radius=0.5*\x,y radius=0.25*\x]; 
        \draw[gray] (-2*\x,\x) arc[start angle=180,end angle=360,x radius=0.5*\x,y radius=0.25*\x]; 
        
        \draw (-1.7*\x,2.4*\x) node[cross=3pt,red] {};
        \draw (-1.3*\x,2.1*\x) node[cross=3pt,red] {};
        \draw (-1.6*\x,1.6*\x) node[cross=3pt,red] {};
        
        \draw (-2*\x, 0) -- (-2*\x, \height);
        \draw (-\x, 0) -- (-\x, \height);
        \draw (-\x, 3.8*\x) node[anchor=west]{$M$};

        \draw [->] (-0.8*\x, 0.5*\height) -- (0.8*\x, 0.5*\height) node[midway,above] {$f$};
        
        \draw (\x, 0) -- (\x, \height);
        \draw (\x, 3.8*\x) node[anchor=west]{$\mathbb{R}$};

        \draw[color=red,fill=red] (\x, 2.4*\x) circle (1pt);
        \draw[color=red,fill=red] (\x, 2.1*\x) circle (1pt);
        \draw[color=red,fill=red] (\x, 1.6*\x) circle (1pt);
    \end{tikzpicture}
\caption{Schematic of a proper Morse function.}
\end{center}
\end{figure}
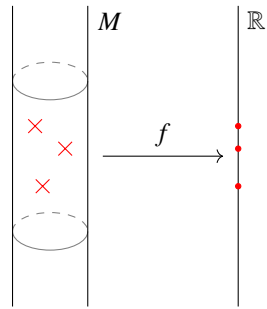

\begin{remark}
There is a subtlety in the presentation of $M$ as a CW complex. If one tries to define the attaching map between cells just by flowing, then the resulting map would be discontinuous. See \cite{M:h-cobo} for the details on the attaching maps.
\end{remark}

The construction of the Morse chain complex for $f$ mirrors the construction of $\Ccell$ of the associated CW complex.
The generators of the Morse complex are the critical points of the Morse function $f$. We let $|x|$ denote the index of the critical point $x$.
For each $x, y \in \Crit f$ we let $\Mod(x,y)$ denote the unparameterised moduli space of flowlines from $x$ to $y$, of dimension $|x|-|y|-1$.
\begin{center}
\resizebox{\textwidth}{!}{
    \begin{tikzpicture}
        \draw[<->] (-2, 0) node[text width=3.5cm,align=center, anchor=east]{Proper Morse function\\ with $|\Crit f|<\infty$}
        -- (2, 0) node[text width=3.5cm,align=center, anchor = west]{Finite based\\ CW complex};
        \draw (0,-5) node[text width=5cm,align=center]{Chain complex\\ with distinguished basis};
        \draw[->] (-3.75,-0.75) 
        to[out=-90, in=135] 
        node[midway, above right]{$\CMorse$} 
        node[midway, above left, text width=3.5cm,align=center, anchor=east]
        {\small\color{red} gen.'s\ $=\Crit f$\\
        \small\color{blue}$\langle \partial x,y \rangle = \#\Mod(x,y)$\\
        \scriptsize\color{orange} (need to orient $\Mod(x,y)$ coherently)} 
        (-2,-4.5);
        \draw[->] (3.75,-0.75) to[out=-90, in=45]
        node[midway, above left]{$\Ccell$} 
        node[midway, above right, text width=3.5cm,align=center, anchor=west]
        {\small\color{red} gen.'s\ = cells\\
        \small\color{blue}$\langle \partial x,y \rangle =$ degree of the attaching map\\
        \scriptsize\color{orange} (need to orient all cells)} 
        (2,-4.5);
    \end{tikzpicture}
}
\end{center}

In the next subsection we will construct framed flow categories, which are `equivalent' to Morse functions with gradient flowlines and to based CW complexes.

\subsection{From Morse functions to framed flow categories}
%
This section is included mostly for motivational purposes, we describe
ideas and constructions which motivate the definitions and results in
Section~\ref{sec:KhovanovHomType}; consequently, no proofs in this
section will be given.

\begin{definition}
    Let $f$ be a Morse-Smale function. We define the \emph{flow category} $\Cat_f$ as the category whose objects are $\Ob\Cat_f = \Crit f$, and whose non-identity morphisms are given by the unparameterised moduli spaces $\Mod(x,y)$ of (possibly broken) flowlines from $x$ to $y$.
\end{definition}

\begin{figure}
\begin{center}
    \begin{tikzpicture}
        \draw (-2, -0.5) -- (-1,0.5) -- (2, 0.5)
        node[anchor=north west]{$\mathbb{D}_x$}
        -- (1, -0.5) -- cycle;
        \draw[<->, color=red] (-.4,-.2) -- (.4, .2);
        \draw[<->] (-.15,.3) -- (.15, -.3);
        \draw[color=black, fill=black] (0,0) node[anchor = south east, pos=20pt]{$x$\,\,} circle (2pt);
        \draw[-{>[sep=2pt]}] (0,0) to[out=-63.5, in=90] (2,-3);
        
        \draw[-] (0,-3) node[anchor=south]{$\mathbb{D}_y$} -- (4,-3);
        \draw[<->, color=red] (2-.45,-3) -- (2+.45, -3);
        \draw[color=black, fill=black] (2,-3) node[anchor=south west]{$y$} circle (2pt);
    \end{tikzpicture}
\end{center}
\caption{An illustration of the equality $T\Mod(x,y) \oplus \R \oplus T\mathbb{D}_y = T\mathbb{D}_x$.}
\label{fig:SARstablenormalbundle}
\end{figure}
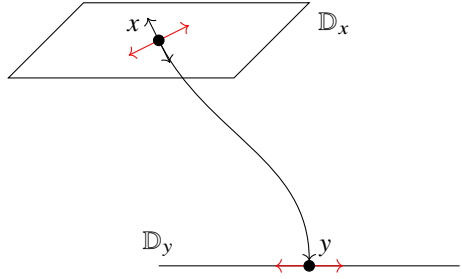

In order to promote $\Cat_f$ to a \emph{framed} flow category we need to endow each $\Mod(x,y)$ with a framing of the stable normal bundle. We do it as follows:
\begin{itemize}
    \item for every downward disc $\mathbb{D}^k$ we choose a trivialisation of $T\mathbb{D}^k$;
    \item for each pair of critical points $x,y$ with $|x|>|y|$ and connected by a flowline, there is an identification
    \[
    T\Mod(x,y) \oplus \R \oplus T\mathbb{D}_y = T\mathbb{D}_x,
    \]
    which gives a stable framing of $T\Mod(x,y)$ (see Figure \ref{fig:SARstablenormalbundle});
    \item a stable framing of $T\Mod(x,y)$ induces a stable framing of the normal bundle too.
\end{itemize}
In the case when $|x|=|y|+1$, the choice of a framing of the stable normal bundle is the same as a (coherent) system of orientations.

\subsection{The hypercube flow category}
This subsection is devoted to a crucial example of flow category, which will be the fundamental building block of the Khovanov stable homotopy type.

The starting point is an innocent-looking function $f \colon \R \to \R$ with only two critical points, at $0$ and $1$, which are respectively a minimum and a maximum. For example, we can choose $f(x) = x^2\cdot\left(\frac32-x\right)$.

\begin{center}
    \begin{tikzpicture}
        \draw[domain=-0.7:1.7,smooth,variable=\x] plot (\x, 1.5*\x*\x-\x*\x*\x);
        \draw[fill=red,color=red] (0,0) circle (2pt);
        \draw[fill=red,color=red] (1,0.5) circle (2pt);
        
        \draw [->] (2, 0.25) -- (3.8, 0.25) node[midway,above] {$f$};
        
        \draw (4, -0.578) -- (4, 1.078);
        \draw (4, 1.1) node[anchor=west]{$\mathbb{R}$};
    \end{tikzpicture}
\end{center}

Define the function $f_n\colon \R^n \to \R$ by setting
\begin{equation}
\label{eq:f_n}
f_n(x_1, \ldots, x_n) = \sum_{i=1}^n f(x_i).
\end{equation}
Then $\Crit f_n = \set{0,1}^n$. The index of a critical point $v \in \set{0,1}^n$ coincide with its $L^1$ norm:
\[
|v| = || v ||_{L^1} = \sum_{i=1}^n v_i.
\]
Moreover, the set $\Crit f_n$ inherits a product partial order by setting $0<1$. For two distinct critical points $u,v \in \set{0,1}^n$, we have that
\[
\Mod(u,v) = \varnothing \qquad \mbox{if $u\not > v$}.
\]
In the next examples we will try to identify $\Mod(u,v)$ when $u>v$.

\begin{exampl3}
\label{ex:function1}
    The function $f=f_1 \colon \R \to \R$ has only two critical points, and the unparameterised moduli space $\Mod(1,0)$ consists of a single point.
\end{exampl3}

\begin{figure}
\begin{center}
    \begin{tikzpicture}
        \draw (0,2) -- (-2,0) -- (0,-2) -- (2,0) -- cycle;

        \foreach \i in {5,...,13}
        {
        \draw[color=cyan, thin] (0,2) to[out=-10*\i, in=10*\i, looseness=1.7] (0,-2);
        }
        
        \draw[fill=black,color=black] (0,2) node[anchor=south]{$11$} circle (2pt);
        \draw[fill=black,color=black] (0,-2) node[anchor=north]{$00$} circle (2pt);
        \draw[fill=blue,color=blue] (-2,0) node[anchor=north east]{$ab$} circle (2pt);
        \draw[fill=blue,color=blue] (2,0) node[anchor=north west]{$cd$} circle (2pt);
        \draw[color=blue, thick] (-2,0) -- (2,0);
        
        \draw (-2,0) node[anchor=south east]{$10$};
        \draw (2,0) node[anchor=south west]{$01$};
        
        
        \draw[fill=red,color=red] (1,1) node[anchor=south west]{$c$} circle (2pt);
        \draw[fill=red,color=red] (1,-1) node[anchor=north west]{$d$} circle (2pt);
        \draw[fill=red,color=red] (-1,1) node[anchor=south east]{$a$} circle (2pt);
        \draw[fill=red,color=red] (-1,-1) node[anchor=north east]{$b$} circle (2pt);
    \end{tikzpicture}
\end{center}
\caption{Illustration of the flowlines from $11$ to $00$ and of the resulting moduli space.}
\label{fig:SARMod1100}
\end{figure}
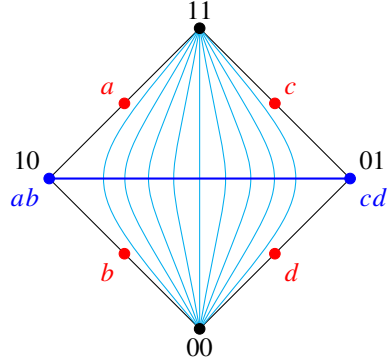

\begin{exampl3}
\label{ex:function2}
    The function $f_2 \colon \R^2 \to \R$ has four critical points, namely $00$, $01$, $10$, and $11$. We have two kinds of moduli spaces, depending on the index difference of $u$ and $v$.

    Using Example \ref{ex:function1}, one can check that if the index difference of two critical points of $f_2$ is exactly 1, then the moduli spaces $\Mod(u,v)$ consists of a single point. We name them
    \[
    \begin{split}
        {\color{red}\Mod(11,10)} = \set{{\color{red}a}}, \qquad
        {\color{red}\Mod(10,00)} = \set{{\color{red}b}}, \\
        {\color{red}\Mod(11,01)} = \set{{\color{red}c}}, \qquad
        {\color{red}\Mod(01,00)} = \set{{\color{red}d}}.
    \end{split}
    \]
    The moduli space $\color{blue}\Mod(11,00)$ is 1-dimensional, and its boundary consists of the two broken flowlines labelled $\color{blue}ab$ and $\color{blue}cd$. The moduli space is in fact an interval with boundary $\set{{\color{blue}ab},{\color{blue}cd}}$.
    One way to visualise $\color{blue}\Mod(11,00)$ is to draw the segment from $10$ to $01$ (in $\R^2$ this segment lies on the line $x+y=1$), and noticing that such a segment intersects each flowline from $11$ to $00$ exactly once. The two endpoints $10$ and $01$ of the segment correspond to the broken flowlines. See also Figure \ref{fig:SARMod1100}.
\end{exampl3}

\begin{figure}
\begin{center}
    \begin{tikzpicture}
    \def\st{1.4}
    
    \def\bl{(-1,-1)}
    \def\br{(1,-1)}
    \def\tl{(-1,1)}
    \def\tr{(1,1)}
    
    \def\bbl{(0,0)}
    \def\bbr{(2,0)}
    \def\btl{(0,2)}
    \def\btr{(2,2)}
    
    \draw \bl node[anchor=north east]{$100$} 
    -- coordinate[midway](a1) \tl node[anchor=south east]{$110$} 
    -- coordinate[midway](a2) \btl node[anchor=south east]{$010$} 
    -- coordinate[midway](a3) \btr node[anchor=south west]{$011$} 
    -- coordinate[midway](a4) \bbr node[anchor=north west]{$001$} 
    -- coordinate[midway](a5) \br node[anchor=north west]{$101$} 
    -- coordinate[midway](a6) \bl;

    \draw[dashed] \bbl node[anchor=south west]{$000$} -- \btl;
    \draw[dashed] \bbl -- \bbr;
    \draw[dashed] \bbl -- \bl;
    
    \foreach \i in {\bl, \tl, \tr, \br, \bbl, \btl, \btr, \bbr}
        {
        \draw[fill=black] \i circle (2pt);
        }

    \draw[fill=purple, fill opacity=0.2, purple] (a1) -- (a2) -- (a3) -- (a4) -- (a5) -- (a6) -- cycle;
    \foreach \i in {1,...,6}
        {
        \draw[purple, fill=purple] (a\i) circle (1pt);
        }

    \draw \tr node[anchor=north west]{$111$} -- \tl;
    \draw \tr -- \br;
    \draw \tr -- \btr;
    
    \end{tikzpicture}
\end{center}
\caption{The moduli space {\color{purple}$\Mod(111,000)$}.}
\label{fig:SARMod111000}
\end{figure}
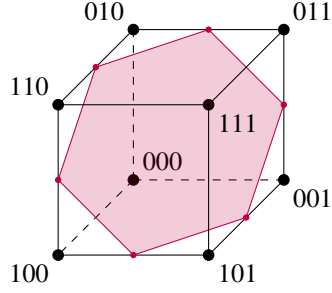

\begin{exercise}
\label{ex:hexagon}
    Identify the flow category associated with $f_3$. Note that $\color{purple}\Mod(111,000)$ is a hexagon: this can be visualised as the section of the unit cube given by the plane $x+y+z=\frac32$, which intersects each (possibly broken) flowline from $111$ to $000$ in exactly one point. See Figure \ref{fig:SARMod111000}.
\end{exercise}

More generally, the moduli space $\Mod(1\cdots1, 0\cdots0)$ in the hypercube flow category is a \emph{permutohedron}.

\begin{definition}
    The $n$-\emph{permutohedron} $\Pi^{n-1}$ is the convex hull of
    \[
    \set{(\sigma(1), \ldots, \sigma(n) ) \in \R^n \,|\, \sigma \in S_n}.
    \]
\end{definition}
Note that for every permutation $\sigma\in S_n$, we have $\sum_{i=1}^n \sigma(i) = \frac{n(n+1)}2$. Thus, $\Pi^{n-1}$ is contained in an $\R^{n-1}$-slice of $\R^n$, and it is in fact an $(n-1)$-dimensional manifold with corners.

The hypercube flow category can also be upgraded to a framed flow category. In order to do so, we need to endow the moduli space $\Mod(u,v)$, which is a permutohedron, with a \emph{coherent framing of the stable normal bundle}.
We hint at its construction by one explicit example.

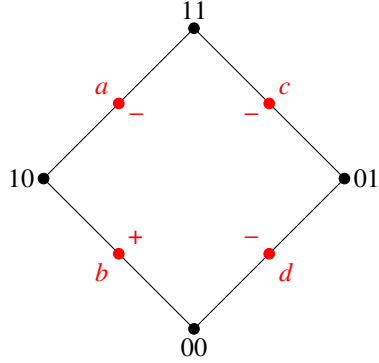
\begin{figure}
\begin{center}
    \begin{tikzpicture}
        \draw (0,2) -- (-2,0) -- (0,-2) -- (2,0) -- cycle;
        
        \draw[fill=black,color=black] (0,2) node[anchor=south]{$11$} circle (2pt);
        \draw[fill=black,color=black] (0,-2) node[anchor=north]{$00$} circle (2pt);
        \draw[fill=black,color=black] (-2,0) node[anchor=east]{$10$} circle (2pt);
        \draw[fill=black,color=black] (2,0) node[anchor=west]{$01$} circle (2pt);
        
        \draw[fill=red,color=red] (1,1)
        node[anchor=south west]{$c$}
        node[anchor=north east]{$-$}
        circle (2pt);
        \draw[fill=red,color=red] (1,-1)
        node[anchor=north west]{$d$}
        node[anchor=south east]{$-$}
        circle (2pt);
        \draw[fill=red,color=red] (-1,1)
        node[anchor=south east]{$a$}
        node[anchor=north west]{$-$}
        circle (2pt);
        \draw[fill=red,color=red] (-1,-1)
        node[anchor=north east]{$b$}
        node[anchor=south west]{$+$}
        circle (2pt);
    \end{tikzpicture}
\end{center}
\caption{A choice of orientation of the 0-dimensional moduli spaces.}
\label{fig:SARorientations}
\end{figure}

Consider the square flow category, i.e.\ the one associated with $f_2$.
A stable framing of the normal bundle of the $0$-dimensional moduli space corresponds to a choice of orientation of each of them (either $+$ or $-$). We need to choose orientations such that exactly one or three of the points $a,b,c,d$ are framed positively; for instance, we may choose orientations as in Figure \ref{fig:SARorientations}.

We envision the moduli spaces $\color{red}\Mod(01,00)$ and $\color{red}\Mod(10,00)$ as sitting in one copy of $\R$, which we denote by $\R_{(i)}$, and the moduli spaces $\color{red}\Mod(11,10)$ and $\color{red}\Mod(11,01)$ as sitting in another copy of $\R$, which we denote by $\R_{(ii)}$.
Thus, we get a product orientation for the broken flowlines $\color{blue}ba$ and $\color{blue}dc$, embedded in $\R_{(i)}\times\R_{(ii)}$, as illustrated in the figure below.

\begin{figure}
\begin{center}
    \begin{tikzpicture}
        \draw[->] (-0.5,0) -- (4,0)
        node[anchor=south] {$\R_{(i)}$};
        \draw[->] (0,-0.5) -- (0,4)
        node[anchor=north west] {$\R_{(ii)}$};

        \foreach \i in {(0,3),(0,1),(1,1),(3,3)}
        {
        \draw[->,cyan] \i -- ++(0,-0.5);
        }
        \foreach \i in {(1,0),(1,1)}
        {
        \draw[->,olive] \i -- ++(0.5,0);
        }
        \foreach \i in {(3,0),(3,3)}
        {
        \draw[->,olive] \i -- ++(-0.5,0);
        }
        
        \draw[fill=red,color=red] (0,1)
        node[anchor=west]{$a$}
        circle (2pt);
        \draw[fill=red,color=red] (0,3)
        node[anchor=west]{$c$}
        circle (2pt);
        \draw[fill=red,color=red] (3,0)
        node[anchor=north]{$d$}
        circle (2pt);
        \draw[fill=red,color=red] (1,0)
        node[anchor=north]{$b$}
        circle (2pt);
        \draw[fill=blue,color=blue] (1,1)
        node[anchor=south]{$ba$}
        circle (2pt);
        \draw[fill=blue,color=blue] (3,3)
        node[anchor=south west]{$dc$}
        circle (2pt);
    \end{tikzpicture}
\end{center}
\caption{The product orientation on the broken flowlines in the square flow category.}
\label{fig:SARproductorientation}
\end{figure}
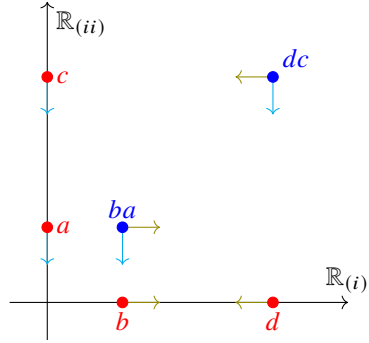

We now wish to orient the stable normal bundle of the moduli space $\color{blue}\Mod(11,00)$, which is an interval with (unoriented) boundary $\set{{\color{blue}ba},{\color{blue}dc}}$. Such a moduli space should come with an embedding
\[
{\color{blue}\Mod(11,00)} \hookrightarrow \R_{(i)} \times \R_{(ii)} \times \R_{\geq0}
\]
whose restriction to the boundary gives the embedding
\[
\set{{\color{blue}ba},{\color{blue}dc}} \hookrightarrow \R_{(i)} \times \R_{(ii)}
\]
constructed above. The coherence condition implies that we can find a framed interval embedded in $\R_{(i)} \times \R_{(ii)} \times \R_{\geq0}$ with boundary the two framed points (see Figure \ref{fig:SARframedinterval}).
The possible stable framings of such an embedded interval are a torsor over
\[
\pi_1(SO):=\varinjlim\,\pi_1(SO(n))\cong\Z/2\Z.
\]

\begin{figure}
\begin{center}
    \begin{tikzpicture}
        \draw[->] (-0.5,0) -- (4,0)
        node[anchor=north] {$\R_{(i)}$};
        \draw[->] (-0.25,-0.25) -- (2,2)
        node[anchor=south east] {$\R_{(ii)}$};
        \draw[->] (0,-0.5) -- (0,4)
        node[anchor=north west] {$\R_{\geq0}$};

        \foreach \i in {(4.5,1.5),(1.5,.5)}
        {
        \draw[->,cyan] \i -- ++(-.25,-0.25);
        }
        \foreach \i in {(4.5,1.5)}
        {
        \draw[->,olive] \i -- ++(0.5,0);
        }
        \foreach \i in {(1.5,.5)}
        {
        \draw[->,olive] \i -- ++(-0.5,0);
        }
        
        \draw[fill=blue,color=blue] (1.5,.5)
        node[anchor=west]{$ba$}
        circle (2pt);
        \draw[fill=blue,color=blue] (4.5,1.5)
        node[anchor=south east]{$dc$}
        circle (2pt);

        \draw[color=blue]
        plot [smooth, tension=1.6]
        coordinates {(4.5,1.5) (3,3) (1.5,.5)};
    \end{tikzpicture}
\end{center}
\caption{An embedding of the moduli space ${\color{blue}\Mod(11,00)}$ into $\R_{(i)} \times \R_{(ii)} \times \R_{\geq0}$.}
\label{fig:SARframedinterval}
\end{figure}

More generally, using obstruction theory we can frame the stable normal bundle coherently for all moduli spaces in the hypercube flow category.

\subsection{The Cohen-Jones-Segal construction}

Cohen-Jones-Segal \cite{CJS} provide a construction to produce a stable based CW complex out of a framed flow category, which we will outline in one particular case in this section.
Their construction shares many similarities with the Pontryagin-Thom construction.

\begin{definition}
    Two topological spaces $X$ and $Y$ are called \emph{stably homotopy equivalent} (denoted $X \sim_S Y$) if $\Sigma^n X$ and $\Sigma^n Y$ are homotopy equivalent for $n \gg 0$.
\end{definition}

\begin{exercise}
\label{ex:T2}
    $\T^2 \sim_S S^1 \vee S^1 \vee S^2$.
    \newline \emph{Hint}: The attaching map of the 2-cell of $\T^2$ is a commutator in $\pi_1$ of the 1-skeleton.
\end{exercise}

\begin{exercise}
    $\CP^2 \not\sim_S S^2 \vee S^4$.
    \newline \emph{Hint}: Steenrod squares are preserved under stable homotopy equivalence.
\end{exercise}

We remark that Exercise \ref{ex:T2} shows that the cup product is not preserved under stable homotopy equivalence.

We will now turn to the construction of a specific example of a stable based CW complex, namely the one associated with the square framed flow category.

Each vertex $v$ of the square will give rise a cell of dimension $|v|+2$.
The resulting stable homotopy type is that of a point.

\begin{center}
    \begin{tikzpicture}
        \draw (0,2) -- (-2,0) -- (0,-2) -- (2,0) -- cycle;
        
        \draw[fill=black,color=black] (0,2) node[anchor=south]{$11$} circle (2pt);
        \draw[fill=black,color=black] (0,-2) node[anchor=north]{$00$} circle (2pt);
        \draw[fill=black,color=black] (-2,0) node[anchor=east]{$10$} circle (2pt);
        \draw[fill=black,color=black] (2,0) node[anchor=west]{$01$} circle (2pt);
        
        \draw[fill=red,color=red] (1,1)
        node[anchor=south west]{$c$}
        node[anchor=north east]{$-$}
        circle (2pt);
        \draw[fill=red,color=red] (1,-1)
        node[anchor=north west]{$d$}
        node[anchor=south east]{$-$}
        circle (2pt);
        \draw[fill=red,color=red] (-1,1)
        node[anchor=south east]{$a$}
        node[anchor=north west]{$-$}
        circle (2pt);
        \draw[fill=red,color=red] (-1,-1)
        node[anchor=north east]{$b$}
        node[anchor=south west]{$+$}
        circle (2pt);

        \begin{scope}[xshift=3.5cm]
            \draw [->] (-2.5,2) -- (0,2) node[anchor=west]{4-cell $C(11)$};
            \draw [->] (-.5,0) -- (0,0) node[anchor=west]{3-cells $C(10)$ and $C(01)$};
            \draw [->] (-2.5,-2) -- (0,-2) node[anchor=west]{2-cell $C(00)$};
        \end{scope}
    \end{tikzpicture}
\end{center}

Given $R \gg \eps > 0$, we parameterise the cells as follows:
\begin{align*}
    C(11) &:= [0,R] \times [-R,R]^2 \times [0,R]; \\
    C(10)=C(01) &:= [0,R] \times [-R,R] \times [-\eps, \eps]; \\
    C(00) &:= [-\eps, \eps]^2.
\end{align*}

\begin{figure}
\begin{center}
    \begin{tikzpicture}
    \def\e{0.5}
    \def\lt{0.25}
    \def\height{3cm}
    
        \draw (0,0) node{$*$};
        \draw[->] (0,-\height+3*\lt cm) --
        node[midway, right] {$\varphi_{00}$}
        (0,-\lt);
        
        \begin{scope}[yshift=-\height]
            \draw[red,fill=red,fill opacity=0.2] 
            (-\e, 0) --
            (0,\e) -- 
            (\e, 0) --
            (0,-\e) -- 
            cycle; 
            \draw[red] (-\e, 0) node[anchor=east]{$C(00)$};
        \end{scope}
    \end{tikzpicture}
\end{center}
\caption{An illustration of the CW complex $X_0$ built as part of the process of constructing the stable homotopy type associated with the square framed flow category.}
\label{fig:SARX0}
\end{figure}
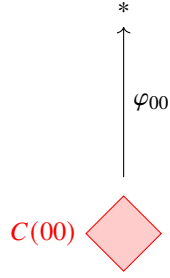

We construct the based CW complex by attaching the cells in order of increasing dimension. We start with the 2-cell $C(00)$, which is attached to the basepoint $*$ via the unique possible attaching map
$\varphi_{00} \colon \partial C(00) \to *$, see Figure \ref{fig:SARX0}.
We denote the resulting CW complex by $X_0$, i.e.
\[
X_0 := C(00) \cup_{\varphi_{00}} \set{*}.
\]

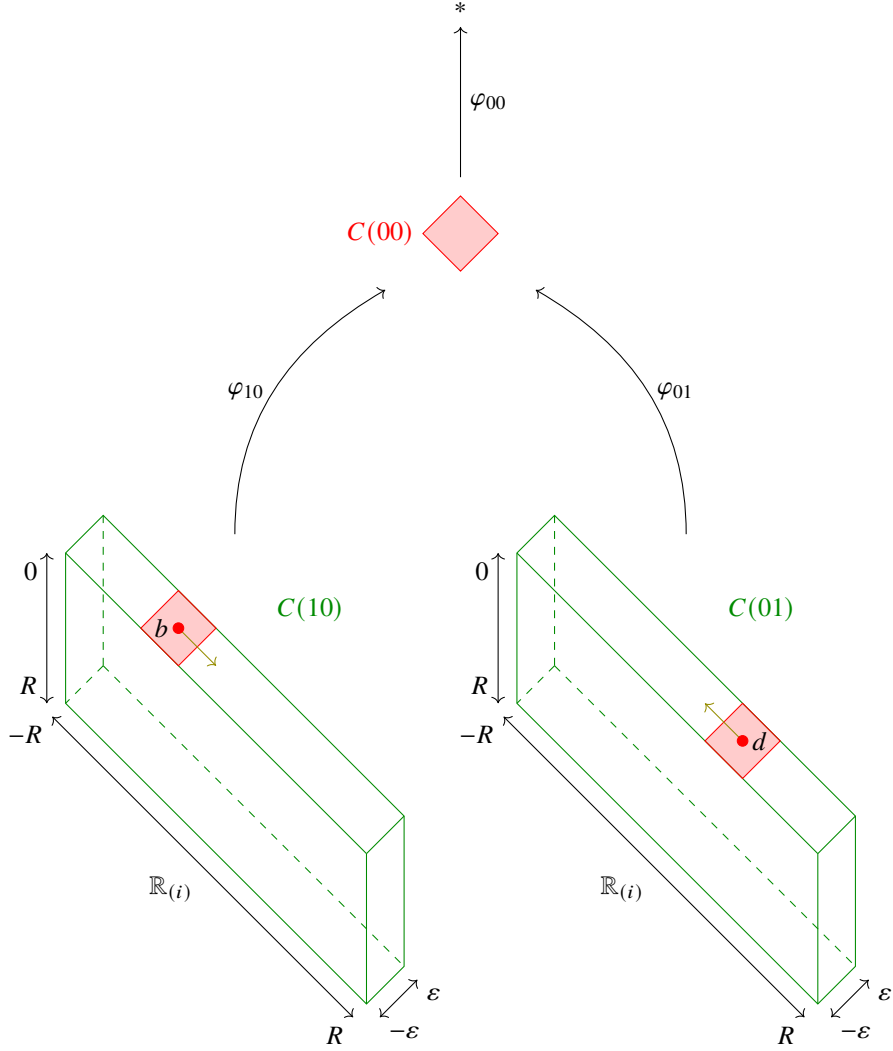
\begin{figure}
\begin{center}
    \begin{tikzpicture}
    \def\srtt{1.41}
    \def\e{0.5}
    \def\he{0.5*\e}
    \def\r{2}
    \def\lt{0.25}
    \def\height{3cm}
    
        \draw (0,0) node{$*$};
        \draw[->] (0,-\height+3*\lt cm) --
        node[midway, right] {$\varphi_{00}$}
        (0,-\lt);
        
        \begin{scope}[yshift=-\height]
            \draw[red,fill=red,fill opacity=0.2] 
            (-\e, 0) --
            (0,\e) -- 
            (\e, 0) --
            (0,-\e) -- 
            cycle; 
            \draw[red] (-\e, 0) node[anchor=east]{$C(00)$};
        \end{scope}

        \foreach \i in {-1,1}
        {
        \begin{scope}[yshift=-3*\height, xshift=\i*\height]

            \draw[Green4] 
            (-\r-\he, -\he) --
            (-\r-\he, \r-\he) --
            (-\r+\he, \r+\he) --
            (\r+\he, -\r+\he) --
            (\r+\he, -2*\r+\he) --
            (\r-\he, -2*\r-\he) -- 
            cycle;
            \draw[Green4] (\r-\he, -\r-\he) -- (-\r-\he, \r-\he);
            \draw[Green4] (\r-\he, -\r-\he) -- (\r+\he, -\r+\he);
            \draw[Green4] (\r-\he, -\r-\he) -- (\r-\he, -2*\r-\he);
            \draw[dashed,Green4] (-\r+\he, \he) -- (-\r-\he, -\he);
            \draw[dashed,Green4] (-\r+\he, \he) -- (-\r+\he, \r+\he);
            \draw[dashed,Green4] (-\r+\he, \he) -- (\r+\he, -2*\r+\he);

            \draw[<->] (-\r-\he-\he, \r-\he) node[anchor=north east]{$0$}
            -- (-\r-\he-\he, -\he) node[anchor=south east]{$R$};
            \draw[<->] (-\r-\he-\he/\srtt, -\he-\he/\srtt) node[anchor=north east]{$-R$}
            -- node[midway, anchor=north east]{$\R_{(i)}$}
            (\r-\he-\he/\srtt, -2*\r-\he-\he/\srtt) node[anchor=north east]{$R$};
            \draw[<->] (\r-\he+\he/\srtt, -2*\r-\he-\he/\srtt) node[anchor=north west]{$-\eps$}
            -- (\r+\he+\he/\srtt, -2*\r+\he-\he/\srtt) node[anchor=north west]{$\eps$};
            
            \begin{scope}[yshift=-1.5*\i*\e cm, xshift=1.5*\i*\e cm]
            \draw[red,fill=red,fill opacity=0.2] 
            (-\e, 0) --
            (0,\e) -- 
            (\e, 0) --
            (0,-\e) -- 
            cycle;

            \draw[->, olive] (0,0) -- (-\i*\e, \i*\e);
            \ifthenelse{\i = -1}
                {\draw[fill=red,color=red] (0,0) node[anchor=east]{\color{black}$b$} circle (2pt);}
                {\draw[fill=red,color=red] (0,0) node[anchor=west]{\color{black}$d$} circle (2pt);}
            \end{scope}

            \ifthenelse{\i = -1}
                {\draw (0.5*\r,0.5*\r) node{\color{Green4}$C(10)$};}
                {\draw (0.5*\r,0.5*\r) node{\color{Green4}$C(01)$};}

            \ifthenelse{\i = -1}
                {
                \draw[->] (0,\r) to[out=90, in=-90+60*\i, looseness=1]
                node[midway, left]{$\varphi_{10}$}
                (-\i*4*\e,2*\height-1.5*\e cm);
                }
                {
                \draw[->] (0,\r) to[out=90, in=-90+60*\i, looseness=1]
                node[midway, right]{$\varphi_{01}$}
                (-\i*4*\e,2*\height-1.5*\e cm);
                }
        \end{scope}
        }
    \end{tikzpicture}
\end{center}
\caption{An illustration of the CW complex $X_0$ built as part of the process of constructing the stable homotopy type associated with the square framed flow category.}
\label{fig:SARX1}
\end{figure}

We now move to the 3-cells. See Figure \ref{fig:SARX1} for reference.
We view the $[-R,R]$ factor as the $\R_{(i)}$ direction, and we recall that we have an embedding of $b$ and $d$ into $\R_{(i)}$, as framed points.

On the boundary of the 3-cell $C(10)$ (resp.\ $C(01)$), we draw a copy of the $[-\eps,-\eps]^2$ on $\set0 \times [-R,R] \times [-\eps,\eps]$, centred at $(0,b,0)$ (resp.\ $(0,d,0)$).

The maps
\begin{align*}
    \varphi_{10} \colon \del C(10) & \to C(00) \cup_{\varphi_{00}} \set{*} \\
    \varphi_{01} \colon \del C(01) & \to C(00) \cup_{\varphi_{00}} \set{*}
\end{align*}
are defined as follows:
\begin{itemize}
    \item on the embedded copies of $[-\eps, \eps]^2$ the maps restrict to degree $\pm1$ maps to $C(00)$, with the sign of the degree depending on the framing of $b$ or $d$;
    \item on the complement of the embedded copy of $[-\eps, \eps]^2$, the maps send everything to the basepoint $\set*$.
\end{itemize}
The maps so defined descend to continuous maps to $C(00) \cup_{\varphi_{00}} \set{*}$.
Using the attaching maps, we construct the CW complex
\[
X_1 := (C(10) \cup C(01) \cup X_0)/\sim.
\]

We now turn to the last remaining cell, namely $C(11)$. Recall that
\[
C(11) = [0,R] \times [-R, R]^2 \times [0,R].
\]
We need to define an attaching map
\[
\varphi_{11} \colon \del C(11) \to X_1.
\]
We view $\partial C(11)$ as the union of three pieces:
\begin{itemize}
    \item {\color{purple} $\set0 \times [-R, R]^2 \times [0,R]$};
    \item {\color{blue} $[0,R] \times [-R, R]^2 \times \set0$};
    \item the part of $\partial C(11)$ with 1st and 4th coordinate both non-zero.
\end{itemize}
Note that the first and the second piece intersect at $\set0 \times [-R,R]^2 \times \set0$.

We consider a copy of $C(01)$ and a copy of $C(10)$ sitting in 
\[
{\color{blue} [0,R] \times [-R, R]^2 \times \set0} \subset \del C(11)
\]
as the slices 
\[
[0,R] \times [-R, R] \times [a-\eps, a+\eps] \times \set0
\qquad
\mbox{and}
\qquad
[0,R] \times [-R, R] \times [c-\eps, c+\eps] \times \set0
\]
respectively, as shown in Figure \ref{fig:SARdelC11blue}.

\begin{figure}
\begin{center}
    \begin{tikzpicture}
    \def\srtt{1.41}
    \def\e{0.5}
    \def\he{0.5*\e}
    \def\r{2}
    \def\lt{0.25}
    \def\height{3cm}

    \foreach \i in {-1,1}
        {
        \begin{scope}[yshift=\i*0.45*\r cm, xshift=\i*0.45*\r cm]

            \draw[Green4] 
            (-\r-\he, -\he) --
            (-\r-\he, \r-\he) --
            (-\r+\he, \r+\he) --
            (\r+\he, -\r+\he) --
            (\r+\he, -2*\r+\he) --
            (\r-\he, -2*\r-\he) -- 
            cycle;
            \draw[Green4] (\r-\he, -\r-\he) -- (-\r-\he, \r-\he);
            \draw[Green4] (\r-\he, -\r-\he) -- (\r+\he, -\r+\he);
            \draw[Green4] (\r-\he, -\r-\he) -- (\r-\he, -2*\r-\he);
            \draw[dashed,Green4] (-\r+\he, \he) -- (-\r-\he, -\he);
            \draw[dashed,Green4] (-\r+\he, \he) -- (-\r+\he, \r+\he);
            \draw[dashed,Green4] (-\r+\he, \he) -- (\r+\he, -2*\r+\he);

            \draw[Green4,fill=Green4,fill opacity=0.3]
            (-\r-\he, \r-\he) --
            (-\r+\he, \r+\he) --
            (\r+\he, -\r+\he) --
            (\r-\he, -\r-\he) --
            cycle;
            \draw[Green4,fill=Green4,fill opacity=0.3]
            (\r+\he, -\r+\he) --
            (\r+\he, -2*\r+\he) --
            (\r-\he, -2*\r-\he) -- 
            (\r-\he, -\r-\he) --
            cycle;
            \draw[Green4,fill=Green4,fill opacity=0.15]
            (\r-\he, -2*\r-\he) --
            (-\r-\he, -\he) --
            (-\r-\he, \r-\he) --
            (\r-\he, -\r-\he) --
            cycle;
            
            \begin{scope}[yshift=-1.5*\i*\e cm, xshift=1.5*\i*\e cm]
            \draw[red,fill=white,fill opacity=1] 
            (-\e, 0) --
            (0,\e) -- 
            (\e, 0) --
            (0,-\e) -- 
            cycle;
            \draw[red,fill=red,fill opacity=0.2] 
            (-\e, 0) --
            (0,\e) -- 
            (\e, 0) --
            (0,-\e) -- 
            cycle;
            
                \draw[->, olive] (0,0) -- (-\i*\e, \i*\e);
                \draw[->, cyan] (0,0) -- (-\e,-\e);
                \ifthenelse{\i = -1}
                {\draw[fill=red,color=red] (0,0) node[anchor=south east]{\color{black}$ba$} circle (2pt);}
                {\draw[fill=red,color=red] (0,0) node[anchor=west]{\color{black}$dc$} circle (2pt);}
            \end{scope}

        \end{scope}
        }

    \begin{scope} 
    
        \draw[blue] 
        (-2*\r,0) --
        (0,2*\r) --
        (2*\r,0) --
        (2*\r,-\r) --
        (0,-3*\r) --
        (-2*\r,-\r) -- 
        cycle;
        
        \draw[blue] (0,-2*\r) -- (-2*\r,0);
        \draw[blue] (0,-2*\r) -- (2*\r,0);
        \draw[blue] (0,-2*\r) -- (0,-3*\r);
        
        \draw[dashed,blue] (0,\r) -- (0,2*\r);
        \draw[dashed,blue] (0,\r) -- (-2*\r,-\r);
        \draw[dashed,blue] (0,\r) -- (2*\r,-\r);

        \draw[<->] (-2*\r-\he,0) node[anchor=north east]{$0$}
        -- (-2*\r-\he, -\r) node[anchor=south east]{$R$};
        \draw[<->] (-2*\r-\he/\srtt, -\r-\he/\srtt) node[anchor=north east]{$-R$}
        -- node[midway, anchor=north east]{$\R_{(i)}$}
        (-\he/\srtt, -3*\r-\he/\srtt) node[anchor=north east]{$R$};
        \draw[<->] (2*\r+\he/\srtt, -\r-\he/\srtt) node[anchor=north west]{$R$}
        -- node[midway, anchor=north west]{$\R_{(ii)}$}
        (\he/\srtt, -3*\r-\he/\srtt) node[anchor=north west]{$-R$};
    
    \end{scope}

        
    \end{tikzpicture}
\end{center}
\caption{The ``blue'' piece ${\color{blue} [0,R] \times [-R, R]^2 \times \set0}$ of $\del C(11)$.}
\label{fig:SARdelC11blue}
\end{figure}

We define the map $\varphi_{11}$ separately on the three pieces of $\del C(11)$.
We start from the blue piece {\color{blue} $[0,R] \times [-R, R]^2 \times \set0$}: with the above picture in mind, it is natural to map the embedded cells $C(01)$ and $C(10)$ inside $\del C(11)$ to the respective abstract ones via degree $\pm1$ maps, depending on the framings of $a$ and $c$,
and to extend this to a map from {\color{blue} $[0,R] \times [-R, R]^2 \times \set0$} to $X_1$ by sending everything in the complement to $*$.

\begin{figure}
\begin{center}
    \begin{tikzpicture}
    \def\srtt{1.41}
    \def\e{0.5}
    \def\he{0.5*\e}
    \def\r{2}
    \def\lt{0.25}
    \def\height{3cm}

        \def\i{-1}
        \begin{scope}[yshift=\i*0.45*\r cm -\r cm, xshift=\i*0.45*\r cm]
        \begin{scope}[yshift=-1.5*\i*\e cm, xshift=1.5*\i*\e cm]
            \draw[red,fill=red,fill opacity=0.2] 
            (-\e, 0) node(a1){} --
            (0,\e) node(d1){} -- 
            (\e, 0) node(c1){} --
            (0,-\e) node(b1){} -- 
            cycle;
            \draw (1.6*\e,1.8*\e) node(c){};
            \draw (0.45*\e,2.6*\e) node(d){};
            
            \draw[->, olive] (0,0) -- (-\i*\e, \i*\e);
            \draw[->, cyan] (0,0) -- (-\e,-\e);
            \ifthenelse{\i = -1}
            {\draw[fill=red,color=red] (0,0) node[anchor=south east]{} circle (2pt);}
            {\draw[fill=red,color=red] (0,0) node[anchor=west]{} circle (2pt);}
        \end{scope}
        \end{scope}

        \def\i{1}
        \begin{scope}[yshift=\i*0.45*\r cm -\r cm, xshift=\i*0.45*\r cm]
        \begin{scope}[yshift=-1.5*\i*\e cm, xshift=1.5*\i*\e cm]
            \draw[red,fill=red,fill opacity=0.2] 
            (-\e, 0) node(b2){} --
            (0,\e) node(c2){} -- 
            (\e, 0) node(d2){} --
            (0,-\e) node(a2){} -- 
            cycle;
            
            \draw[->, olive] (0,0) -- (-\i*\e, \i*\e);
            \draw[->, cyan] (0,0) -- (-\e,-\e);
            \ifthenelse{\i = -1}
            {\draw[fill=red,color=red] (0,0) node[anchor=south east]{} circle (2pt);}
            {\draw[fill=red,color=red] (0,0) node[anchor=west]{} circle (2pt);}
        \end{scope}
        \end{scope}

        \draw[-] (a1.center) to[out=90, in=90, looseness=1.5] (a2.center);
        \draw[-] (b1.center) to[out=90, in=90, looseness=1.5] (b2.center);
        \draw[-] (c1.center) to[out=90, in=55, looseness=1.2] (c.center);
        \draw[dashed,-] (c.center) to[out=55, in=90, looseness=1.2] (c2.center);
        \draw[dashed,-] (d1.center) to[out=90, in=60, looseness=1.2] (d.center);
        \draw[-] (d.center) to[out=60, in=90, looseness=1.2] (d2.center);

    \begin{scope} 
    
        \draw[purple] 
        (-2*\r,0) --
        (0,2*\r) --
        (2*\r,0) --
        (2*\r,-\r) --
        (0,-3*\r) --
        (-2*\r,-\r) -- 
        cycle;
        
        \draw[purple] (0,-2*\r) -- (-2*\r,0);
        \draw[purple] (0,-2*\r) -- (2*\r,0);
        \draw[purple] (0,-2*\r) -- (0,-3*\r);
        
        \draw[dashed,purple] (0,\r) -- (0,2*\r);
        \draw[dashed,purple] (0,\r) -- (-2*\r,-\r);
        \draw[dashed,purple] (0,\r) -- (2*\r,-\r);

        \draw[<->] (-2*\r-\he,0) node[anchor=north east]{$R$}
        -- (-2*\r-\he, -\r) node[anchor=south east]{$0$};
        \draw[<->] (-2*\r-\he/\srtt, -\r-\he/\srtt) node[anchor=north east]{$-R$}
        -- node[midway, anchor=north east]{$\R_{(i)}$}
        (-\he/\srtt, -3*\r-\he/\srtt) node[anchor=north east]{$R$};
        \draw[<->] (2*\r+\he/\srtt, -\r-\he/\srtt) node[anchor=north west]{$R$}
        -- node[midway, anchor=north west]{$\R_{(ii)}$}
        (\he/\srtt, -3*\r-\he/\srtt) node[anchor=north west]{$-R$};
    
    \end{scope}
    
        
    \end{tikzpicture}
\end{center}
\caption{The ``purple'' piece ${\color{purple} \set0 \times [-R, R]^2 \times [0,R]}$ of $\del C(11)$.}
\label{fig:SARdelC11purple}
\end{figure}
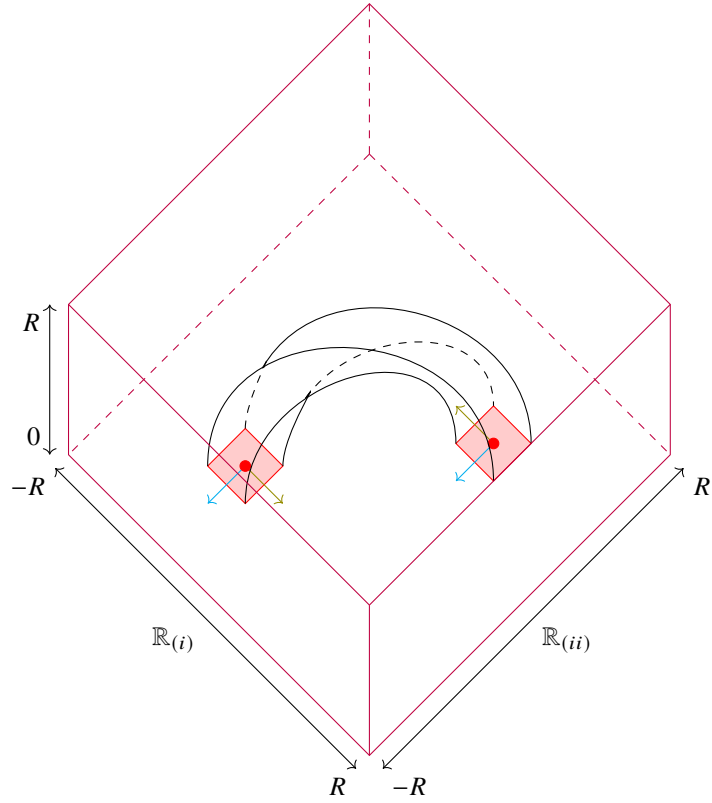

Then, a na\"{i}ve approach to extend $\varphi_{11}$ to all of $\del
C(11)$ would be to map the purple and the black pieces to the
basepoint $*$. However, the map so defined would not be continuous,
because the whole purple piece of $\del C(11)$ would be mapped to $*$,
except the two red squares in $\set0 \times [-R,R]^2 \times \set0$,
which lie in $\del({\color{purple}\set0 \times [-R, R]^2 \times
  [0,R]})$.  Thus, we modify the map in the purple portion of $\del
C(11)$ by using the framed interval from $ab$ to $cd$. We envision the
framed interval constructed above as sitting inside
${\color{purple}\set0 \times [-R, R]^2 \times [0,R]}$, where the last
coordinate corresponds to $\R_{\geq0}$. A regular neighbourhood of it
can be parameterised by $[-\eps,\eps]^2 \times I$.  We can then define
the map $\varphi_{11}$ on ${\color{purple}\set0 \times [-R, R]^2
  \times [0,R]}$ by setting it on the neighbourhood of the framed
interval as
\[
[-\eps,\eps]^2 \times I \xrightarrow{pr_1} [-\eps,\eps]^2 \xrightarrow{\sim} C(00) \rightarrow X_0 \subset X_1,
\]
and by sending the complement to $*$. See Figure \ref{fig:SARdelC11purple}.

Thus, we have so far defined $\varphi_{11}$ on {\color{purple} $\set0 \times [-R, R]^2 \times [0,R]$} and {\color{blue} $[0,R] \times [-R, R]^2 \times \set0$}. We extend it to $\del C(11)$ by sending the rest of $\del C(11)$ to $*$.

\section{A Khovanov stable homotopy type}
\label{sec:KhovanovHomType}

Our next goal is to construct a Khovanov homotopy type starting from the Khovanov complex. We quickly recall the construction of the Khovanov complex. For a more detailed exposition we refer the reader to \cite{K:Khovanov, BN:Khovanov}.

\subsection{The Khovanov complex}

\begin{figure}
    \centering
\begingroup%
  \makeatletter%
  \providecommand\color[2][]{%
    \errmessage{(Inkscape) Color is used for the text in Inkscape, but the package 'color.sty' is not loaded}%
    \renewcommand\color[2][]{}%
  }%
  \providecommand\transparent[1]{%
    \errmessage{(Inkscape) Transparency is used (non-zero) for the text in Inkscape, but the package 'transparent.sty' is not loaded}%
    \renewcommand\transparent[1]{}%
  }%
  \providecommand\rotatebox[2]{#2}%
  \newcommand*\fsize{\dimexpr\f@size pt\relax}%
  \newcommand*\lineheight[1]{\fontsize{\fsize}{#1\fsize}\selectfont}%
  \ifx\svgwidth\undefined%
    \setlength{\unitlength}{255.11794987bp}%
    \ifx\svgscale\undefined%
      \relax%
    \else%
      \setlength{\unitlength}{\unitlength * \real{\svgscale}}%
    \fi%
  \else%
    \setlength{\unitlength}{\svgwidth}%
  \fi%
  \global\let\svgwidth\undefined%
  \global\let\svgscale\undefined%
  \makeatother%
  \begin{picture}(1,0.11111101)%
    \lineheight{1}%
    \setlength\tabcolsep{0pt}%
    \put(0,0){\includegraphics[width=\unitlength,page=1]{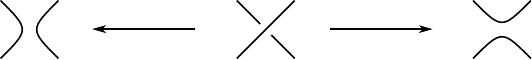}}%
    \put(0.2673437,0.07062701){\color[rgb]{0,0,0}\makebox(0,0)[lt]{\lineheight{1.25}\smash{\begin{tabular}[t]{l}$0$\end{tabular}}}}%
    \put(0.70243709,0.07062701){\color[rgb]{0,0,0}\makebox(0,0)[lt]{\lineheight{1.25}\smash{\begin{tabular}[t]{l}$1$\end{tabular}}}}%
  \end{picture}%
\endgroup%

    \caption{The 0-smoothing and the 1-smoothing of a crossing.}
    \label{fig:smoothings}
\end{figure}

\begin{figure}
    \centering
\begingroup%
  \makeatletter%
  \providecommand\color[2][]{%
    \errmessage{(Inkscape) Color is used for the text in Inkscape, but the package 'color.sty' is not loaded}%
    \renewcommand\color[2][]{}%
  }%
  \providecommand\transparent[1]{%
    \errmessage{(Inkscape) Transparency is used (non-zero) for the text in Inkscape, but the package 'transparent.sty' is not loaded}%
    \renewcommand\transparent[1]{}%
  }%
  \providecommand\rotatebox[2]{#2}%
  \newcommand*\fsize{\dimexpr\f@size pt\relax}%
  \newcommand*\lineheight[1]{\fontsize{\fsize}{#1\fsize}\selectfont}%
  \ifx\svgwidth\undefined%
    \setlength{\unitlength}{45.50228761bp}%
    \ifx\svgscale\undefined%
      \relax%
    \else%
      \setlength{\unitlength}{\unitlength * \real{\svgscale}}%
    \fi%
  \else%
    \setlength{\unitlength}{\svgwidth}%
  \fi%
  \global\let\svgwidth\undefined%
  \global\let\svgscale\undefined%
  \makeatother%
  \begin{picture}(1,0.98819127)%
    \lineheight{1}%
    \setlength\tabcolsep{0pt}%
    \put(0,0){\includegraphics[width=\unitlength,page=1]{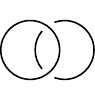}}%
    \put(0.435,0.81803755){\color[rgb]{0,0,0}\makebox(0,0)[lt]{\lineheight{1.25}\smash{\begin{tabular}[t]{l}$1$\end{tabular}}}}%
    \put(0.435,-0.05){\color[rgb]{0,0,0}\makebox(0,0)[lt]{\lineheight{1.25}\smash{\begin{tabular}[t]{l}$2$\end{tabular}}}}%
  \end{picture}%
\endgroup%

    \caption{A non-minimal diagram of the 2-component unlink, together with a labelling of the crossings.}
    \label{fig:unlink_diag}
\end{figure}

\begin{figure}
    \centering
    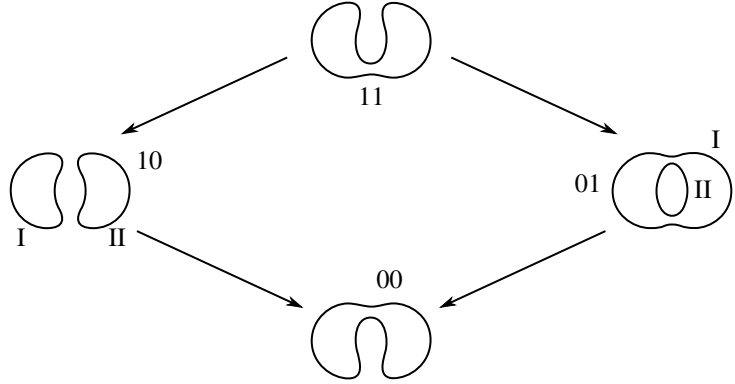
    \caption{The cube of resolutions of the diagram for the unlink in Figure \ref{fig:unlink_diag}.
    We have also fixed a total ordering on the circles at each resolution: the resolutions, indicated by the labels $\mathrm{I}$ and $\mathrm{II}$ for the circles at the resolutions $10$ and $01$. (The other resolutions have only one circle, so there is only one total order there.)}
    \label{fig:unlink_cube}
  \end{figure}

Given a link diagram and a crossing, we define the 0-smoothing and the 1-smoothing of the crossing as in Figure \ref{fig:smoothings}.
Let $D$ be a link diagram with $n$ crossings, labelled in some order from $1$ to $n$. See Figure \ref{fig:unlink_diag} for an example.
Each $v \in \set{0,1}^n$ determines a \emph{complete resolution} $D_v$ by smoothing crossing $i$ according to the $i$-th entry of $v$.
Topologically, each complete resolution consists of an unlink. We can organise the complete resolutions on a hypercube of dimension $n$, so that each $v$ corresponds to a vertex of the hypercube.
See Figure \ref{fig:unlink_cube}.

Let the \emph{weight} or \emph{index} $|v|$ of $v$ be the number of coordinates of $v$ with value $1$. Two vertices $u$ and $v$ are connected by an edge $e$ of the hypercube if and only if $u$ and $v$ differ only at one entry, and in such a case the edge is oriented from the vertex with \emph{higher} weight to the vertex with \emph{lower} weight. Note that this convention is opposite to the standard one in the literature (e.g., \cite{K:Khovanov, BN:Khovanov}); equivalently, in these notes we consider Khovanov homology as opposed to the more popular Khovanov cohomology.

From the cube of resolution Khovanov constructed a complex as follows. Let $V= \Z\langle 1,x\rangle$ be the $\Z$-module spanned by two elements $1$ and $x$. To each $v \in \set{0,1}^n$ associate the vector space
\[
V_v = \bigotimes_{\mathrm{components\,of\,} D_v} V,
\]
where each factor $V$ is considered as associated to a link component of $D_v$.
The module underlying the Khovanov complex is defined as
\[
\CKh (D) = \bigoplus_{v \in \set{0,1}} V_{v}.
\]
Recall that in the construction of framed flow category \emph{spacifying} a given chain complex $C$ one needs to specify a basis of $C$ as well, see Remark \ref{rem:spac-basis}.
Luckily for us, the $\Z$-module $\CKh(D)$ comes with a preferred basis. For each $v \in \set{0,1}$, the $\Z$-module $V_{v}$ has a basis consisting of all possible tensor products of the generators $1$ and $x$ for $V$; we call them the \emph{standard generators} of $V_v$, and we denote them with a $v$ subscript, e.g. if $D_v$ has four components, after totally ordering them, the element $1x11_{v}$ is the Khovanov generator that labels the second component by $x$ and other components by $1$.
As $v$ varies in $\set{0,1}$, all standard generators give a basis for $\CKh(D)$.

In order to define the differential on $\CKh(D)$, Khovanov defined two maps $m$ and $\Delta$, associated to the merge and split cobordisms respectively:
\begin{align*}
    m \colon V \otimes V & \to V &
    \Delta \colon V & \to V \otimes V \\
    1 \otimes 1 & \mapsto 1 &
    1 & \mapsto 1 \otimes x + x \otimes 1 \\
    1 \otimes x & \mapsto x &
    x & \mapsto x \otimes x \\
    x \otimes 1 & \mapsto x \\
    x \otimes x & \mapsto 0     
\end{align*}

If there is an edge from $u$ to $v$, then there is an elementary cobordism between $D_u$ and $D_v$, which is a split saddle or a merge saddle. We define a map
\[
\del_{u,v} \colon V_u \to V_v
\]
as the map $m$ or $\Delta$ on the $V$ factors associated with the components of $D_u$ and $D_v$ on which the saddle cobordism happens, and as the identity map on all other $V$ factors.

Fix a framing of all the 0-dimensional moduli spaces in the hypercube flow category associated to the function $f_n$ defined in Equation \eqref{eq:f_n} at page \pageref{eq:f_n}, such that on every square face, exactly one or three of the points are framed positively.
The differential is then defined by summing over all edges:
\[
\del:=\sum_{\substack{|u-v|=1\\|u|=|v|+1}} (-1)^{s(u,v)} \del_{u,v},
\]
where $s(u,v)=0$ if $\Mod(u,v)$ is framed positively and $s(u,v)=1$ if $\Mod(u,v)$ is framed negatively.

\begin{figure}
    \centering
    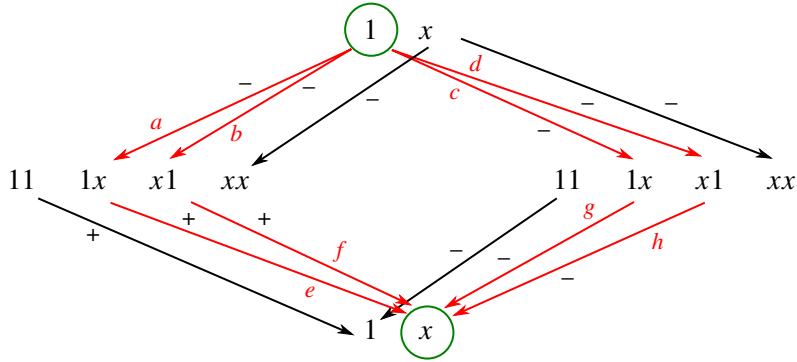
    \caption{
    The Khovanov complex associated with the diagram for the 2-component unlink in Figure \ref{fig:unlink_diag}.
    To simplify the notation, we suppressed the tensor product symbol $\otimes$ and the subscripts specifying the vertex of the hypercube.
    The relative location of the standard generators indicates the vertex of the hypercube to which it belongs. Compare with Figure \ref{fig:unlink_cube}. For example, the four standard generators located at the left are in $V_{(10)}$; among them, the generator $1x$ is the one that labels the first circle (marked I) by $1$ and the second circle (marked II) by $x$.}
    \label{fig:unlink_cx}
\end{figure}

As an illustration, Figure \ref{fig:unlink_cx} shows the Khovanov complex for the diagram of the 2-component unlink in Figure \ref{fig:unlink_diag}. The figure shows the 12 standard generators of the Khovanov complex, organised to mirror the structure of the cube of resolutions (Figure \ref{fig:unlink_cube}). The arrows encode the differential: if there is an arrow from $y$ to $z$, then $z$ appears in the expansion of $\del y$ with coefficient $\pm1$ depending on the sign label on the arrow. (For the moment disregard the fact that some arrows are coloured in red and some of the generators are circled in green.)
The homology of $(\CKh(D), \del)$ is the Khovanov homology of the link. However, we do not take homology, but instead construct a framed flow category from the complex.

\subsection{The Lipshitz-Sarkar framed flow category}

We will use the example from Figures \ref{fig:unlink_diag} and \ref{fig:unlink_cube} to describe the framed flow category $\mc X(D)$ from \cite{LS:KhovanovHomotopyType}.

We start from the set of objects, $\Ob\mc X(D)$, which consists of all standard Khovanov generators.
Thus, for our recurring example of the diagram of the 2-component unlink in Figure \ref{fig:unlink_diag}, the \textbf{objects} are the 12 generators that are shown in Figure \ref{fig:unlink_cx}.

The \textbf{0-dimensional moduli spaces} are framed points, given by the Khovanov differential. Following the convention in Figure \ref{fig:unlink_cx}, the elements of 0-dimensional moduli spaces are the arrows shown in the figure, together with a framing specified by the sign label: for example,
\[
\Mod(x_{(11)}, xx_{(01)}) = \set{p_-}, \qquad 
\Mod(x_{(11)}, x1_{(01)}) = \varnothing,
\]
where $p_-$ is a negatively framed point.
More precisely, given two objects / standard generators $\y$ and $\z$ with $|\y|=|\z|+1$, we represent $\del \y \in CKh(D)$ as a linear combination of standard generators, and we take the coefficient $\langle \del \y, \z \rangle$ of $\z$ in such an expansion. Note that by the definition of the Khovanov differential $\del$, the only values $\langle \del \y, \z \rangle$ can take are $+1$, $-1$, and $0$. Then, we set
\[
\Mod(\y,\z)=
\begin{cases}
    \set{\mbox{a positively framed point}} & \mbox{if $\langle \del \y, \z \rangle = +1$}\\
    \set{\mbox{a negatively framed point}} & \mbox{if $\langle \del \y, \z \rangle = -1$}\\
    \varnothing & \mbox{if $\langle \del \y, \z \rangle = 0$}
\end{cases}
\]

We will construct the higher moduli spaces inductively as disjoint unions of permutohedra, using the fact that $\del \Mod$ is the union of products of lower-dimensional moduli spaces already constructed.

Most of the higher moduli spaces are constructed without making choices, but there is one interesting case, which appears during the construction of the \textbf{1-dimensional moduli spaces}, where a choice is needed. This is the case of $\Mod(1_{(11)}, x_{(00)})$ illustrated in Figure \ref{fig:unlink_cx}: $1_{(11)}$ and $x_{(00)}$ are the two standard generators circled in green; the broken flowlines between them, which constitute the boundary of the moduli space, are obtained by following the red arrows. Thus, we can see that
\[
\del \Mod(1_{(11)}, x_{(00)}) = \set{ea, fb, gc, hd},
\]
where $a$, $b$, $c$, $d$, $e$, $f$, $g$, and $h$ are the red arrows in Figure \ref{fig:unlink_cx}. Of the four broken flowlines, two (namely $ea$ and $fb$) factor through the vertex $10$, and the other two (namely $gc$ and $hd$) factor through the vertex $01$. Note that the vertex through which the broken flowline factors determines its framing.

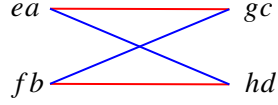
\begin{figure}
    \centering
\begingroup%
  \makeatletter%
  \providecommand\color[2][]{%
    \errmessage{(Inkscape) Color is used for the text in Inkscape, but the package 'color.sty' is not loaded}%
    \renewcommand\color[2][]{}%
  }%
  \providecommand\transparent[1]{%
    \errmessage{(Inkscape) Transparency is used (non-zero) for the text in Inkscape, but the package 'transparent.sty' is not loaded}%
    \renewcommand\transparent[1]{}%
  }%
  \providecommand\rotatebox[2]{#2}%
  \newcommand*\fsize{\dimexpr\f@size pt\relax}%
  \newcommand*\lineheight[1]{\fontsize{\fsize}{#1\fsize}\selectfont}%
  \ifx\svgwidth\undefined%
    \setlength{\unitlength}{99.72479272bp}%
    \ifx\svgscale\undefined%
      \relax%
    \else%
      \setlength{\unitlength}{\unitlength * \real{\svgscale}}%
    \fi%
  \else%
    \setlength{\unitlength}{\svgwidth}%
  \fi%
  \global\let\svgwidth\undefined%
  \global\let\svgscale\undefined%
  \makeatother%
  \begin{picture}(1,0.32597913)%
    \lineheight{1}%
    \setlength\tabcolsep{0pt}%
    \put(0,0){\includegraphics[width=\unitlength,page=1]{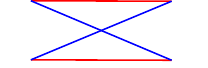}}%
    \put(-0.00074188,0.29290223){\color[rgb]{0,0,0}\makebox(0,0)[lt]{\lineheight{1.25}\smash{\begin{tabular}[t]{l}\large$ea$\end{tabular}}}}%
    \put(-0.00074174,0.00865532){\color[rgb]{0,0,0}\makebox(0,0)[lt]{\lineheight{1.25}\smash{\begin{tabular}[t]{l}\large$fb$\end{tabular}}}}%
    \put(0.88208169,0.29290223){\color[rgb]{0,0,0}\makebox(0,0)[lt]{\lineheight{1.25}\smash{\begin{tabular}[t]{l}\large$gc$\end{tabular}}}}%
    \put(0.88208169,0.00865541){\color[rgb]{0,0,0}\makebox(0,0)[lt]{\lineheight{1.25}\smash{\begin{tabular}[t]{l}\large$hd$\end{tabular}}}}%
  \end{picture}%
\endgroup%

    \caption{The two possibilities for $\Mod(1_{(11)}, x_{(00)})$ from the example in Figure \ref{fig:unlink_cx}. The {\color{red} \emph{right} choice} is shown in red, and the {\color{blue} \emph{left} choice} is shown in blue.
    The two possible matchings are called \emph{ladybug} matchings, and a choice between the two of them must be done each time a ladybug configuration as in Figure \ref{fig:ladybug_diag} appears.}
    \label{fig:ladybug_choice}
\end{figure}

Thus, if we want it to be a union of permutohedra, $\Mod(1_{(11)}, x_{(00)})$ should consist of two intervals, which can be chosen in two sensible ways: either an interval connecting $ea$ to $gc$ and one connecting $fb$ to $hd$ (the {\color{red}\emph{right} choice}), or an interval connecting $ea$ to $hd$ and one connecting $fb$ to $gc$ (the {\color{blue}\emph{left} choice}). See Figure \ref{fig:ladybug_choice}.
In a bombastic moment of self-confidence, we claim that \underline{we always make the right choice}!
The 1-dimensional moduli spaces so defined can be upgraded to \emph{framed} moduli spaces, using a framing on the cube flow category.

\begin{figure}
    \centering
    \includegraphics{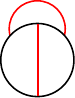}
    \caption{A ladybug matching. The figure represents an all-1 resolution of a link diagram (in fact, of the diagram in Figure \ref{fig:unlink_diag}). The red arcs, one for each crossing of the link diagram, trace the band surgeries that turn each 1-smoothing into a 0-smoothing.}
    \label{fig:ladybug_diag}
\end{figure}

\begin{remark}
In the language of \cite{LS:KhovanovHomotopyType}, we just made a choice of a \emph{ladybug matching}. The name comes from the fact that the diagram $D$ can be alternatively specified by its all-$1$ resolution, together with a collection of embedded arcs that record the band surgeries needed to turn a 1-smoothing into a 0-smoothing.
Figure \ref{fig:ladybug_diag} shows the all-1 resolution and the surgery arcs for the diagram in Figure \ref{fig:unlink_diag}. The shape resembling a ladybug motivates the name \emph{ladybug configuration}.

Each time a ladybug configuration appears we must choose a ladybug matching. The fact that we always make the right choice will be important for the next step, namely the construction of the 2-dimensional moduli spaces.
\end{remark}

Surprisingly, we will see that the choice of the ladybug matchings on the 1-dimensional moduli spaces is \emph{the only choice we need} to make, under the assumption that each moduli space $\Mod(\y,\z)$ is a disjoint union of permutohedra.

\begin{figure}
    \centering
    \includegraphics{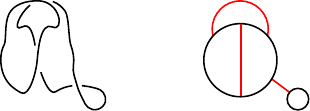}
    \caption{The figure on the left shows a link diagram for the 2-component unlink. The figure on the right represents its all-1 resolution together with the arcs that turn each 1-smoothing into a 0-smoothing.}
    \label{fig:headphones_diag}
\end{figure}

\begin{figure}
    \centering
    \includegraphics{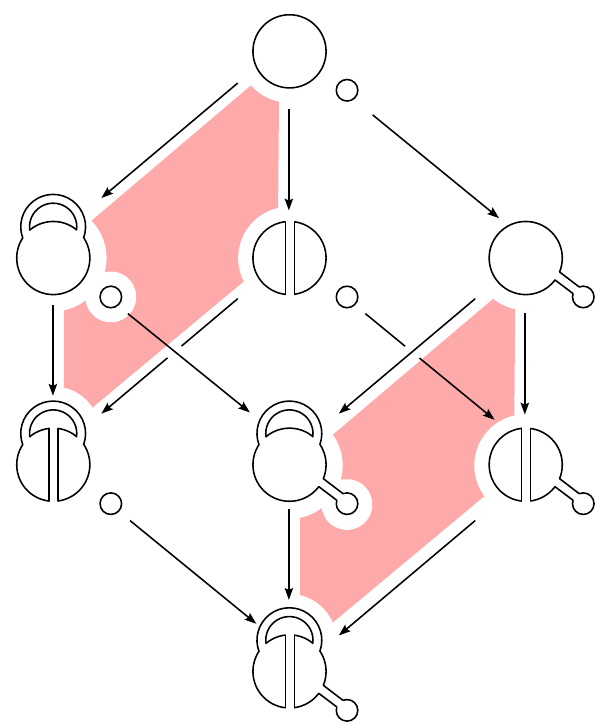}
    \caption{The cube of resolutions for the link diagram in Figure \ref{fig:headphones_diag}. The faces shaded in red correspond to ladybug configurations.}
    \label{fig:headphones_cube}
\end{figure}

We now turn to the construction of the \textbf{2-dimensional moduli spaces}. We start from the example in Figure \ref{fig:headphones_diag}, whose associated cube of resolutions is shown in Figure \ref{fig:headphones_cube}. One can check that the 3-dimensional cube of resolutions has 2 ladybug faces, which are shaded in light red.

\begin{figure}
    \centering
    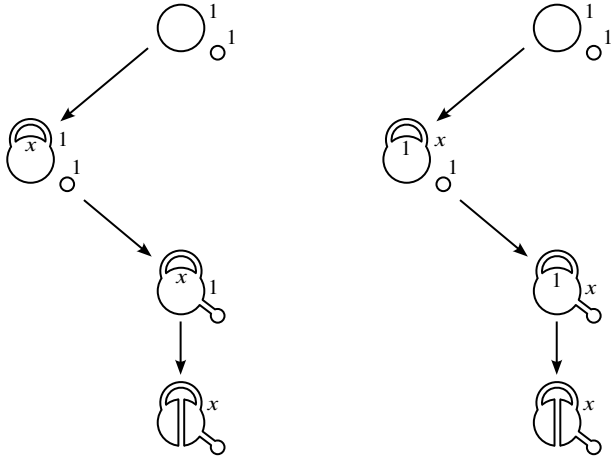
    \caption{The figure shows two fully broken flowlines in the boundary of a 2-dimensional moduli space arising from the cube of resolutions associated to the link diagram in Figure \ref{fig:headphones_diag}.}
    \label{fig:headphones_fullybroken}
\end{figure}

We consider the moduli space
\begin{center}
\begingroup%
  \makeatletter%
  \providecommand\color[2][]{%
    \errmessage{(Inkscape) Color is used for the text in Inkscape, but the package 'color.sty' is not loaded}%
    \renewcommand\color[2][]{}%
  }%
  \providecommand\transparent[1]{%
    \errmessage{(Inkscape) Transparency is used (non-zero) for the text in Inkscape, but the package 'transparent.sty' is not loaded}%
    \renewcommand\transparent[1]{}%
  }%
  \providecommand\rotatebox[2]{#2}%
  \newcommand*\fsize{\dimexpr\f@size pt\relax}%
  \newcommand*\lineheight[1]{\fontsize{\fsize}{#1\fsize}\selectfont}%
  \ifx\svgwidth\undefined%
    \setlength{\unitlength}{98.36235827bp}%
    \ifx\svgscale\undefined%
      \relax%
    \else%
      \setlength{\unitlength}{\unitlength * \real{\svgscale}}%
    \fi%
  \else%
    \setlength{\unitlength}{\svgwidth}%
  \fi%
  \global\let\svgwidth\undefined%
  \global\let\svgscale\undefined%
  \makeatother%
  \begin{picture}(1,0.28416969)%
    \lineheight{1}%
    \setlength\tabcolsep{0pt}%
    \put(0,0){\includegraphics[width=\unitlength,page=1]{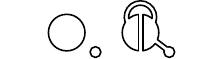}}%
    \put(-0.00174673,0.07915239){\color[rgb]{0,0,0}\makebox(0,0)[lt]{\lineheight{1.25}\smash{\begin{tabular}[t]{l}\large$\mathcal{M}($\end{tabular}}}}%
    \put(0.91434947,0.0791523){\color[rgb]{0,0,0}\makebox(0,0)[lt]{\lineheight{1.25}\smash{\begin{tabular}[t]{l}\large$).$\end{tabular}}}}%
    \put(0.39099207,0.21543921){\color[rgb]{0,0,0}\makebox(0,0)[lt]{\lineheight{1.25}\smash{\begin{tabular}[t]{l}\footnotesize$1$\end{tabular}}}}%
    \put(0.44137389,0.07875559){\color[rgb]{0,0,0}\makebox(0,0)[lt]{\lineheight{1.25}\smash{\begin{tabular}[t]{l}\footnotesize$1$\end{tabular}}}}%
    \put(0.79317337,0.14496522){\color[rgb]{0,0,0}\makebox(0,0)[lt]{\lineheight{1.25}\smash{\begin{tabular}[t]{l}\footnotesize$x$\end{tabular}}}}%
    \put(0.51674261,0.06390266){\color[rgb]{0,0,0}\makebox(0,0)[lt]{\lineheight{1.25}\smash{\begin{tabular}[t]{l}\large{,}\end{tabular}}}}%
  \end{picture}%
\endgroup%

\end{center}
Its boundary (already defined, since it is the union of products of lower-dimensional moduli spaces) contains 12 fully broken flowlines:
\begin{itemize}
    \item 4 fully broken flowlines for each ladybug face;
    \item 2 fully broken flowlines factoring through
    \begin{center}
\begingroup%
  \makeatletter%
  \providecommand\color[2][]{%
    \errmessage{(Inkscape) Color is used for the text in Inkscape, but the package 'color.sty' is not loaded}%
    \renewcommand\color[2][]{}%
  }%
  \providecommand\transparent[1]{%
    \errmessage{(Inkscape) Transparency is used (non-zero) for the text in Inkscape, but the package 'transparent.sty' is not loaded}%
    \renewcommand\transparent[1]{}%
  }%
  \providecommand\rotatebox[2]{#2}%
  \newcommand*\fsize{\dimexpr\f@size pt\relax}%
  \newcommand*\lineheight[1]{\fontsize{\fsize}{#1\fsize}\selectfont}%
  \ifx\svgwidth\undefined%
    \setlength{\unitlength}{90.40145514bp}%
    \ifx\svgscale\undefined%
      \relax%
    \else%
      \setlength{\unitlength}{\unitlength * \real{\svgscale}}%
    \fi%
  \else%
    \setlength{\unitlength}{\svgwidth}%
  \fi%
  \global\let\svgwidth\undefined%
  \global\let\svgscale\undefined%
  \makeatother%
  \begin{picture}(1,0.30919414)%
    \lineheight{1}%
    \setlength\tabcolsep{0pt}%
    \put(0,0){\includegraphics[width=\unitlength,page=1]{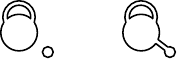}}%
    \put(0.37046214,0.11125198){\color[rgb]{0,0,0}\makebox(0,0)[lt]{\lineheight{1.25}\smash{\begin{tabular}[t]{l}and\end{tabular}}}}%
    \put(0.98438841,0.11125198){\color[rgb]{0,0,0}\makebox(0,0)[lt]{\lineheight{1.25}\smash{\begin{tabular}[t]{l},\end{tabular}}}}%
  \end{picture}%
\endgroup%

    \end{center}
    namely the ones illustrated in Figure \ref{fig:headphones_fullybroken};
    \item 2 fully broken flowlines factoring through
    \begin{center}
\begingroup%
  \makeatletter%
  \providecommand\color[2][]{%
    \errmessage{(Inkscape) Color is used for the text in Inkscape, but the package 'color.sty' is not loaded}%
    \renewcommand\color[2][]{}%
  }%
  \providecommand\transparent[1]{%
    \errmessage{(Inkscape) Transparency is used (non-zero) for the text in Inkscape, but the package 'transparent.sty' is not loaded}%
    \renewcommand\transparent[1]{}%
  }%
  \providecommand\rotatebox[2]{#2}%
  \newcommand*\fsize{\dimexpr\f@size pt\relax}%
  \newcommand*\lineheight[1]{\fontsize{\fsize}{#1\fsize}\selectfont}%
  \ifx\svgwidth\undefined%
    \setlength{\unitlength}{90.40145514bp}%
    \ifx\svgscale\undefined%
      \relax%
    \else%
      \setlength{\unitlength}{\unitlength * \real{\svgscale}}%
    \fi%
  \else%
    \setlength{\unitlength}{\svgwidth}%
  \fi%
  \global\let\svgwidth\undefined%
  \global\let\svgscale\undefined%
  \makeatother%
  \begin{picture}(1,0.23898171)%
    \lineheight{1}%
    \setlength\tabcolsep{0pt}%
    \put(0,0){\includegraphics[width=\unitlength,page=1]{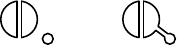}}%
    \put(0.37046214,0.11125186){\color[rgb]{0,0,0}\makebox(0,0)[lt]{\lineheight{1.25}\smash{\begin{tabular}[t]{l}and\end{tabular}}}}%
    \put(0.98438841,0.11125186){\color[rgb]{0,0,0}\makebox(0,0)[lt]{\lineheight{1.25}\smash{\begin{tabular}[t]{l},\end{tabular}}}}%
  \end{picture}%
\endgroup%

    \end{center}
    which the reader is encouraged to find by themselves.
\end{itemize}

The boundary of such a moduli space is therefore a 1-dimensional manifold with 12 corners (one for each fully broken trajectory). It turns out that it must be either a dodecagon (i.e., 12-gon) or the disjoint union of two hexagons (i.e., 6-gon). Lipshitz and Sarkar \cite{LS:KhovanovHomotopyType} showed that if we choose all the ladybug matchings to be right (or all of them to be left), then we get two hexagons. Thus, we can define the moduli space
\begin{center}

\end{center}
as the disjoint union of two filled hexagons, which are permutohedra (see Exercise \ref{ex:hexagon}). Once again, these inherit framings from a framing of the cube flow category.

Finally, we turn to the \textbf{moduli spaces of dimension $\geq 3$}. Their definition is actually easier than in the cases of lower dimension.

Assume that $\y$ and $\z$ are standard generators of the Khovanov chain complex, such that $|\y| > |\z| + 3$.
By induction, we assume that the boundary of the moduli space $\Mod(\y,\z)$ is already defined. It turns out that $\del \Mod(\y,\z)$ is a regular cover of the moduli space in the cube flow category, i.e.\ a permutohedron.
If $n := |\y| - |\z| > 3$, then the permutohedron $\Pi^{n-1}$ has simply connected boundary, and therefore every regular cover of it is trivial. Thus, there is no obstruction to defining the moduli space $\Mod(\y,\z)$ as a disjoint union of permutohedra of dimension $n-1$, which also inherit framings from a framing of the cube flow category.

\section{A knot Floer stable homotopy type}

Once we have the `spacification'
of Khovanov homology, and 
witness its success, it is
natural to ask a similar construction for invariants
of knots defined in the Heegaard Floer context. As the
construction requires a good understanding of `moduli spaces',
it is instructive to start the process
with the version of knot Floer homology where the 
most control on the moduli spaces is available: grid homology.
For an introduction to grid homology see
\cite{MOS, MOSzT, OSS}.
For general (pointed or doubly pointed) Heegaard
diagrams, the determination of the moduli spaces required for
describing the flow category is analytically extremely challenging, and indeed in general these spaces depend on further choices
(such as an almost complex structure, or perturbations), 
which are hard to incorporate into the theory. 
These difficulties can be bypassed by the combinatorial 
description provided by grid homology.

Despite the simplicity of the description of the chain complex
involved in the definition of grid homology, the presence of `bubbles'
provides serious difficulties for setting up the stable homotopy
type. These bubbles do arise already on the homological level, but
they cancel out, hence in grid homology they are `invisible' to the
boundary map, resulting in a (relatively) simple theory. In the
construction of the framed flow category, however, the bubbles cannot
be ignored anymore, bringing in an extra level of complexity in the
theory.  In this section (following \cite{ManolescuSarkar}) the
construction of a spectrum with stable homology giving grid homology
is outlined. The appropriate topological invariance of the spectrum,
that is, the proof of the fact that the result is a knot-link
invariant is, however, still open. This invariance would provide the
real power of the theory, since it would open the way to construct
further knot invariants, which might go beyond the mere homology
theory given by grid homology.

\subsection{Grid diagrams}
Recall that every knot (and even link) in the standard three-sphere
$S^3$ can be presented by a grid diagram, a combinatorial object  with the following properties:

Consider the grid of $n^2$ small squares in an $n\times n $ 
square in the plane for some $n\in {\mathbb {N}}$. Suppose that each column and each row has exactly one marking $X$ and exactly one marking
$O$ in one of the small squares in the column or row at hand.  By connecting the $X$ in a column to the $O$ in the same
column, and the $O$ in the row to the $X$ in the same row 
(oriented by this order), and with the convention that a
vertical segment always passes over a horizontal one, the grid diagram
provides the diagram of a (PL embedded) link in $S^3$. It is not hard to see
that every knot and link can be presented in this way, and 
a method similar to the proof of Reidemeister's fundamental
result in knot theory provides two types of  moves with the
property that grids defining isotopic knots/links can 
be transformed into each other by a sequence of these `grid moves'.

Indeed, by identifying the top and bottom edges of the 
big square, and doing the same with the left-most and
right-most edges, we get a grid diagram on the standard 
torus of $S^3$. The two grid moves of commutation and 
stabilization then provide a similar theory
for knots in this context.

\subsection{Grid homology}

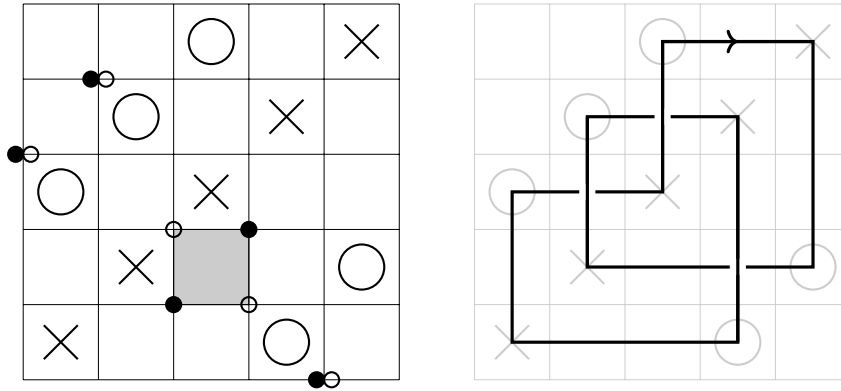
\begin{figure}[t]
\centering
\begin{tikzpicture}
\begin{scope}[xshift=-3cm]
\draw[step=1cm,black,very thin] (0,0) grid (5,5);
%
\drawO{1}{3}
\drawO{2}{4}
\drawO{3}{5}
\drawO{4}{1}
\drawO{5}{2}
%
\drawX{1}{1}
\drawX{2}{2}
\drawX{3}{3}
\drawX{4}{4}
\drawX{5}{5}
%

\draw[very thin, fill=black, fill opacity=0.2] (2,1) rectangle (3,2);

\draw[thick, fill] (-0.1,3) circle (0.1);
\draw[thick, fill] (0.9,4) circle (0.1);
\draw[thick, fill] (2,1) circle (0.1);
\draw[thick, fill] (3,2) circle (0.1);
\draw[thick, fill] (3.9,0) circle (0.1);
\draw[thick] (0.1,3) circle (0.1);
\draw[thick] (1.1,4) circle (0.1);
\draw[thick] (2,2) circle (0.1);
\draw[thick] (3,1) circle (0.1);
\draw[thick] (4.1,0) circle (0.1);

\end{scope}

\begin{scope}[xshift=3cm]
  \begin{scope}[color=black, opacity=0.2]
\draw[step=1cm,very thin] (0,0) grid (5,5);
%
\drawO{1}{3}
\drawO{2}{4}
\drawO{3}{5}
\drawO{4}{1}
\drawO{5}{2}
%
\drawX{1}{1}
\drawX{2}{2}
\drawX{3}{3}
\drawX{4}{4}
\drawX{5}{5}
\end{scope}

\draw[very thick]
(0.5,0.5) --
(3.5,0.5) --
(3.5,3.5) --
(1.5,3.5) --
(1.5,1.5) --
(4.5,1.5) --
(4.5,4.5) --
(2.5,4.5) --
(2.5,2.5) --
(0.5,2.5) --
cycle;

\begin{scope}[very thick,decoration={
    markings,
    mark=at position 0.5 with {\arrow{>}}}
    ] 
    \draw[postaction={decorate}] (2.5,4.5) -- (4.5,4.5);
\end{scope}


\draw[fill=white, white] (3.4,1.4) rectangle (3.6,1.6);
\draw[fill=white, white] (1.4,2.4) rectangle (1.6,2.6);
\draw[fill=white, white] (2.4,3.4) rectangle (2.6,3.6);


\draw[very thick] (1.5,1.5) -- (1.5,3.5);
\draw[very thick] (2.5,2.5) -- (2.5,4.5);
\draw[very thick] (3.5,0.5) -- (3.5,3.5);
\end{scope}
\end{tikzpicture}
\caption{The left hand side shows a grid diagram for the left-handed trefoil knot illustrated on the right. The full and the hollow dots determine two grid states $\x$ and $\y$ respectively. A small square (out of the 25 ones for this grid) is shaded. The shaded small square is also a rectangle from $\x$ to $\y$.}
\end{figure}

The chain complex of grid homology is generated by 
grid states: a \emph{grid state} $\x$ is simply a 
bijection between the horizontal and vertical circles
in the toroidal grid. A geometric presentation of such
an object can be given by choosing an intersection point
of the vertical and horizontal circles in the grid in 
such a way that each vertical and each horizontal circle
contains exactly one chosen point. The differential in 
the chain complex is given by counting ``empty rectangles'',
cf. Figure~\ref{fig:egy}.
In more detail, to a pair of two distinct grid states $\x,\y$ we
associate 0 if $\x$ and $\y$ differ in more than 2 coordinates.
If the difference is in exactly two circles, then the
circles of these coordinates (there are two vertical and two horizontal of those)
partition the torus into four rectangles, and by choosing
an appropriate orientation convention, two of them will point
from $\x$ to $\y$ and two from $\y$ to $\x$. Now we count the
rectangles from $\x$ to $\y$ which are \emph{empty} (i.e., which do not contain
any further coordinate of $\x$) and 
do not contain $X$-marking. The actual set-up of the
homology theory is slightly more involved (requires a system
of weights determined by the $O$-markings, and consequently
is defined over a more elaborate ring); we will not go into the
details of the theory here, and advise the interested reader
to consult \cite{OSS} for a more complete account.
We just mention here that a crucial step in building 
the theory is to show that the boundary map 
$\partial$ (counting 
empty rectangles with some weights) has square zero, 
and therefore gives rise to a chain complex, and hence to
homologies. The geometric input in the proof of
$\partial ^2=0$ is an incarnation of Gromov's compactness
result, showing that the compositions of two rectangles
(counted by $\partial ^2$) come in pairs. (Over the
field of two elements this fact is sufficient, over the
integers a sign assignment should be also fixed so these
pairs contribute zero to the final count.)
In turn, this pairing follows from a simple planar geometric
observation decomposing the union of two rectangles in
two different ways --- this is depicted in 
Figures~\ref{fig:ketto} and \ref{fig:harom}. The situation
is somewhat more complicated in further versions of the theory
(where we allow the empty rectangle to contain $X$-markings), 
as
the concatenation of two rectangles might add up to a 
strip around the torus; this contribution will be 
cancelled by another such strip, as shown in 
Figure~\ref{fig:tizenegy}. The significance of this last 
small problem becomes more crucial in the spacification 
process, as it will be hinted later on.

The spacification of grid homology will once again be done by
constructing a framed flow category $\Cat$ associated with a grid
diagram. The objects of the category are the grid states, i.e.\ the
$n$-tuples of intersection points of the vertical and horizontal
circles such that each vertical and each horizontal circle contains
exactly one chosen point.  Given two grid states $\x$ and $\y$, the
moduli space $\Mod(\x,\y)$ will be a smooth manifold with corners. By
contrast with the Khovanov homotopy type, the construction is more
involved and the moduli spaces will not always be disjoint unions of
permutohedra. In the most general case, the moduli space $\Mod(\x,\y)$
will be obtained by gluing together the moduli spaces $\Mod(D)$
associated with all the \emph{domains} from $\x$ to $\y$ with Maslov
index $|\x|-|\y|$.  Recall from \cite{OSS} that a domain from $\x$ to
$\y$ in a given grid diagram is a formal linear combination
\[
D=\sum n_{i,j} S_{i,j}
\]
of
(the closures) of the small squares $\{ S_{i,j}\}_{i,j=1, \ldots , n}$
with the following properties:
\begin{itemize}
\item the horizontal part of
the boundary (i.e. the part of the boundary $\partial D$ as a 2-chain on the
horizontal circles) is a 1-chain $\partial _hD$ with
$\partial (\partial _hD)=\y -\x$, while
\item the vertical part of the boundary $\partial _vD$ satisfies
  $\partial (\partial _vD)=\x -\y$.
\end{itemize}
The domain $D=\sum n_{i,j} S_{i,j}$ is \emph{positive} if all $n_{i,j}\geq 0$.

The Maslov index $\mu (D)$ of a domain $D$ is the formal dimension of
the space of holomorphic representatives of the domain (when
considered as a 2-chain); it admits a rather simple formula, which
will be sufficient in our subsequent discussion.  For a domain
$D=\sum n_{i,j}S_{i,j}$ from $\x $ to $\y$ we define the point measure
$p(D)$ as follows: for each $x_i\in \x$ (i.e.\ an intersection point of
the horizontal and vertical lines which belongs to $\x$) we
take the average of the four numbers $n_{i,j}$ corresponding to the four
small squares meeting at $x_i$.  Then $p_{\x}(D)$ is the sum of these
local contributions for all $x_i\in \x$, and $p(D)=p_{\x}(D)+p_{\y}(D)$. The
Maslov index of $D$ is then
\[
\mu (D)=p(D).
\]
\begin{remark}
  The general formula for the Maslov index in Heegaard Floer theory
  also involves a term originating from the geometry of $D$ (called
  the Euler measure of $D$); as this term is additive and vanishes for
  the small squares, in a grid diagram all domains have vanishing
  Euler measure.
  \end{remark}
  
In the next subsections we will give examples of the moduli spaces
associated with some domains, and in Subsection~\ref{ssec:GeneralCase}
we describe an inductive way (based on the Maslov index) for
constructing these spaces.

\subsection{Maslov index 1 domains}
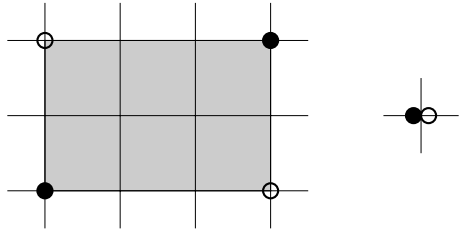
\begin{figure}[h]
\centering
\begin{tikzpicture}

\begin{scope}[xshift=-3cm]
\draw[step=1cm,black,very thin] (-0.5,-0.5) grid (3.5,2.5);
\draw[very thin, fill=black, fill opacity=0.2] (0, 0) rectangle (3,2);
\draw[thick, fill] (0,0) circle (0.1);
\draw[thick, fill] (3,2) circle (0.1);
\draw[thick] (0,2) circle (0.1);
\draw[thick] (3,0) circle (0.1);
\end{scope}

\begin{scope}[xshift=2cm, yshift=1cm]
\draw[step=1cm,black,very thin] (-0.5,-0.5) grid (0.5,0.5);
\draw[thick, fill] (-0.1,0) circle (0.1);
\draw[thick] (0.1,0) circle (0.1);
\end{scope}

\end{tikzpicture}
\caption{A rectangle connecting the generators $\x$ (full circle) and $\y$ (hollow circle). Note that the 
two generators have all further coordinates
equal, symbolized by the crossing on the right.}
\label{fig:egy}
\end{figure}

These are the domains we need to encounter in the
boundary map of grid homology. 
When $\x$ and $\y$ can be connected by a rectangle (that is,
they differ in two circles),
the moduli space
$\Mod (\x, \y)$
contains two, one or zero points
depending on the number of empty rectangles connecting $\x$
and $\y$. An example of an empty rectangle is shown by 
Figure~\ref{fig:egy}.

This case corresponds to the connecting domains with 
Maslov index $\mu (D)=1$.
The space associated to such a domain is a single point, hence
the moduli space $\Mod (\x, \y)$ (the morphisms between the two
objects $\x$ and $\y$ in the flow category) 
consists of two, one or zero points.

\subsection{Maslov index 2 domains}
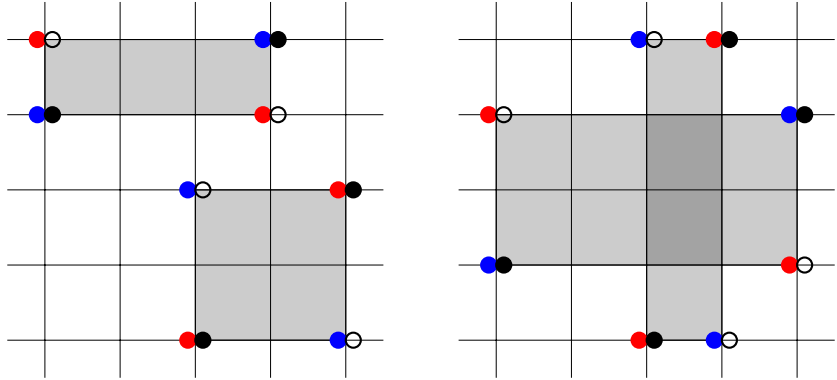
\begin{figure}[h]
\centering
\begin{tikzpicture}

\begin{scope}[xshift=-3cm]
\draw[step=1cm,black,very thin] (-0.5,-0.5) grid (4.5,4.5);
\draw[very thin, fill=black, fill opacity=0.2] (0, 3) rectangle (3,4);
\draw[very thin, fill=black, fill opacity=0.2] (2, 0) rectangle (4,2);
%
\draw[thick, red, fill] (0-0.1,4) circle (0.1);
\draw[thick, red, fill] (3-0.1,3) circle (0.1);
\draw[thick, red, fill] (2-0.1,0) circle (0.1);
\draw[thick, red, fill] (4-0.1,2) circle (0.1);
\draw[thick, blue, fill] (0-0.1,3) circle (0.1);
\draw[thick, blue, fill] (3-0.1,4) circle (0.1);
\draw[thick, blue, fill] (2-0.1,2) circle (0.1);
\draw[thick, blue, fill] (4-0.1,0) circle (0.1);
\draw[thick, fill] (0+0.1,3) circle (0.1);
\draw[thick, fill] (3+0.1,4) circle (0.1);
\draw[thick, fill] (2+0.1,0) circle (0.1);
\draw[thick, fill] (4+0.1,2) circle (0.1);
\draw[thick] (0+0.1,4) circle (0.1);
\draw[thick] (3+0.1,3) circle (0.1);
\draw[thick] (2+0.1,2) circle (0.1);
\draw[thick] (4+0.1,0) circle (0.1);
\end{scope}

\begin{scope}[xshift=3cm]
\draw[step=1cm,black,very thin] (-0.5,-0.5) grid (4.5,4.5);
\draw[very thin, fill=black, fill opacity=0.2] (0, 1) rectangle (4,3);
\draw[very thin, fill=black, fill opacity=0.2] (2, 0) rectangle (3,4);
\draw[thick, red, fill] (0-0.1,3) circle (0.1);
\draw[thick, red, fill] (4-0.1,1) circle (0.1);
\draw[thick, red, fill] (2-0.1,0) circle (0.1);
\draw[thick, red, fill] (3-0.1,4) circle (0.1);
\draw[thick, blue, fill] (0-0.1,1) circle (0.1);
\draw[thick, blue, fill] (4-0.1,3) circle (0.1);
\draw[thick, blue, fill] (2-0.1,4) circle (0.1);
\draw[thick, blue, fill] (3-0.1,0) circle (0.1);
\draw[thick, fill] (0+0.1,1) circle (0.1);
\draw[thick, fill] (4+0.1,3) circle (0.1);
\draw[thick, fill] (2+0.1,0) circle (0.1);
\draw[thick, fill] (3+0.1,4) circle (0.1);
\draw[thick] (0+0.1,3) circle (0.1);
\draw[thick] (4+0.1,1) circle (0.1);
\draw[thick] (2+0.1,4) circle (0.1);
\draw[thick] (3+0.1,0) circle (0.1);
\end{scope}

\end{tikzpicture}
\caption{
Each of the two figures show a Maslov index 2 domain from the (black) full circle to the hollow circle. Geometrically, the domain consists of two rectangles, which can be disjoint (as in the left) or overlapping (as on the right).
In either case, the Maslov index 2 domain can be expressed as the concatenation of two rectangles (i.e., doing one rectangle after the other) in two different ways, depending on which rectangle comes first. One such concatenation factors through the grid state given by red circles, and the other one factors through the grid state given by blue circles. The resulting moduli space is a segment, as in Figure \ref{fig:negy}, with the endpoints (broken flowlines) given by the two concatenations just described.}
\label{fig:ketto}
\end{figure}

 In the next step we consider 
domains with Maslov index $\mu (D)=2$ connecting the grid states $\x$ and $\y$. These domains are 
relatively easy to understand: they are concatenations of 
two rectangles. A simple analysis shows that this can happen in  two ways. In the first case the two rectangles `move' different coordinate pairs (so the two rectangles move four coordinates 
in total). Geometrically in this case the domain $D$ is the
union of two rectangles (which might not be disjoint, but they
do not share sides), see Figure~\ref{fig:ketto}.

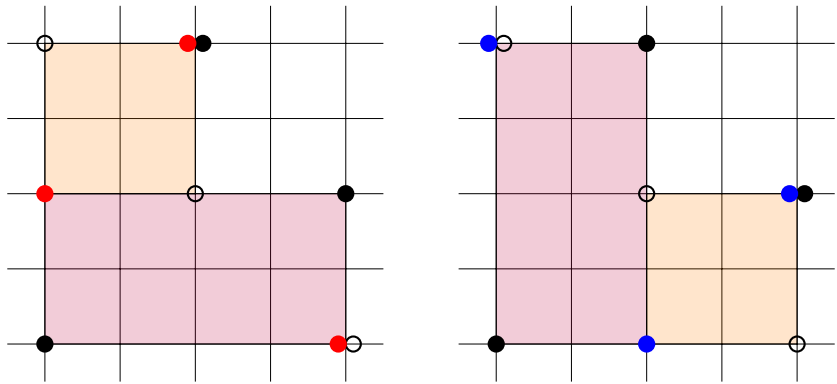
\begin{figure}[h]
\centering
\begin{tikzpicture}

\begin{scope}[xshift=-3cm]
\draw[step=1cm,black,very thin] (-0.5,-0.5) grid (4.5,4.5);
\draw[very thin, fill=purple, fill opacity=0.2] (0, 0) rectangle (4,2);
\draw[very thin, fill=orange, fill opacity=0.2] (0, 2) rectangle (2,4);
\draw[thick, fill] (0,0) circle (0.1);
\draw[thick, fill] (4,2) circle (0.1);
\draw[thick, fill] (2+0.1,4) circle (0.1);
\draw[thick] (0,4) circle (0.1);
\draw[thick] (4.1,0) circle (0.1);
\draw[thick] (2,2) circle (0.1);
\draw[thick, red, fill] (2-0.1,4) circle (0.1);
\draw[thick, red, fill] (4-0.1,0) circle (0.1);
\draw[thick, red, fill] (0,2) circle (0.1);
\end{scope}

\begin{scope}[xshift=3cm]
\draw[step=1cm,black,very thin] (-0.5,-0.5) grid (4.5,4.5);
\draw[very thin, fill=purple, fill opacity=0.2] (0, 0) rectangle (2,4);
\draw[very thin, fill=orange, fill opacity=0.2] (2, 0) rectangle (4,2);
\draw[thick, fill] (0,0) circle (0.1);
\draw[thick, fill] (4.1,2) circle (0.1);
\draw[thick, fill] (2,4) circle (0.1);
\draw[thick] (0.1,4) circle (0.1);
\draw[thick] (4,0) circle (0.1);
\draw[thick] (2,2) circle (0.1);
\draw[thick, blue, fill] (-0.1,4) circle (0.1);
\draw[thick, blue, fill] (4-0.1,2) circle (0.1);
\draw[thick, blue, fill] (2,0) circle (0.1);
\end{scope}

\end{tikzpicture}
\caption{The $L$-shaped region decomposes into two rectangles in two different ways.
Such a region is a domain from $\x$ (full circles) to $\y$ (hollow circles); the two decompositions factor through $\w$ (red circles) or $\z$ (blue circles).}
\label{fig:harom}
\end{figure}

In the second case one of the moving coordinates of the 
first rectangle is also a moving coordinate of the second one.
In this case, in geometric terms the domain $D$ is an L-shaped
region in the grid torus, see Figure~\ref{fig:harom}.

The domain $D$ can be decomposed into the concatenation of two empty
rectangles in two different ways, once as a rectangle $D_1$ from $\x$
to a grid state $\w$ composed with a rectangle $D_2$ from $\w$ to
$\y$, and also as the composition of $D_1'$ from $\x$ to $\z$ with
$D_2'$ from $\z$ to $\y$. When the two rectangles move four
coordinates, the geometric presentation is quite obvious: $D_1$ (as a
subset of the grid torus) agrees with $D_2'$ and $D_2$ with
$D_1'$. When the rectangles share moving coordinates, the
decomposition is given by the two ways an L-shaped domain decomposes
as the union of two rectangles.  This is indicated in
Figure~\ref{fig:harom}.

\begin{figure}[h]
\centering
\begin{tikzpicture}

\begin{scope}[xshift=-2cm]
\node (x) at (0, 2) {$\x$};
\node (y) at (0,-2) {$\y$};
\node (z) at (-1,0) {$\z$};
\node (w) at (1,0) {$\w$};
\draw[->, thick] (x) -- (z) node [midway, left] {$a$};
\draw[->, thick] (x) -- (w) node [midway, right] {$c$};
\draw[->, thick] (z) -- (y) node [midway, left] {$b$};
\draw[->, thick] (w) -- (y) node [midway, right] {$d$};
\end{scope}

\begin{scope}[xshift=3cm]
\draw (-2,0) node {$\mathcal{M}(D) \,\, =$};
\draw[thick, fill] (-1,0) node (X) {} circle (0.1);
\draw (-1,-0.1) node[below] {$ab$};
\draw[thick, fill] (1,0) node (Y) {} circle (0.1);
\draw (1,-0.1) node[below] {$cd$};
\draw[thick] (X) -- (Y);
\end{scope}

\end{tikzpicture}
\caption{The figure on the left shows the broken flowlines for a domain $D$ of Maslov index 2 (such as those in Figures \ref{fig:ketto} and \ref{fig:harom}).
The figure on the right show the moduli space $\Mod(D)$.}
\label{fig:negy}
\end{figure}
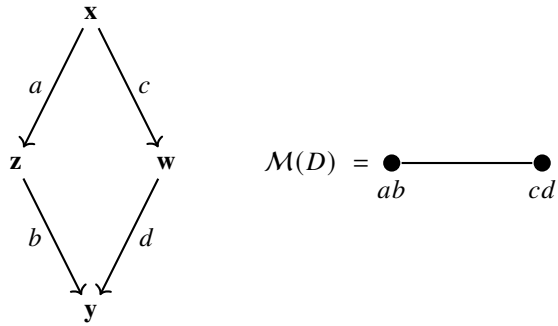

Indeed, these decompositions correspond to the broken
flowline picture, hence we associate to such a domain
an interval as moduli space, with the two endpoints corresponding to the 
two products of one-point moduli spaces, as it is shown in 
Figure~\ref{fig:negy}. The moduli space in this case associated to
the domain $D$ is a manifold with boundary (an interval, really), and
$\Mod (\x , \y )$ will be (as always) the union of all the spaces
associated to Maslov index 2 domains from $\x$ to $\y$.

\subsection{Maslov index 3 domains}
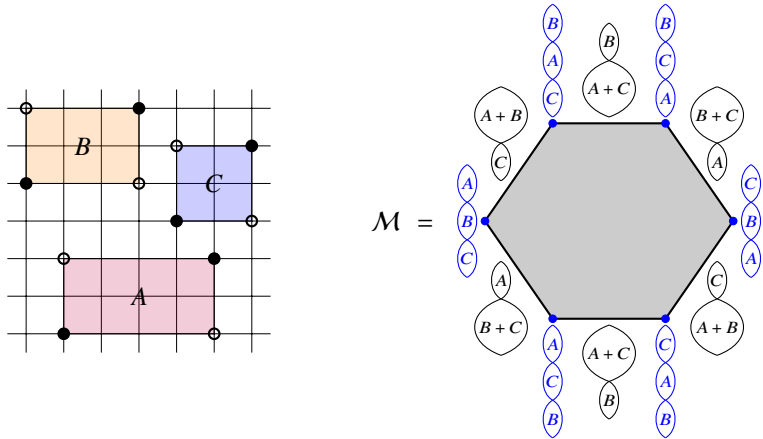
\begin{figure}[h]
\centering
\begin{tikzpicture}

\begin{scope}[xshift=-3cm, yshift=-1.5cm]
\draw[step=0.5cm,black,very thin] (-0.25,-0.25) grid (3.25,3.25);
\draw[very thin, fill=purple, fill opacity=0.2] (0.5, 0) rectangle (2.5,1) node[midway, opacity=1](C-red){$A$};
\draw[very thin, fill=orange, fill opacity=0.2] (0, 2) rectangle (1.5,3) node[midway, opacity=1](C-orange){$B$};
\draw[very thin, fill=blue, fill opacity=0.2] (2, 1.5) rectangle (3,2.5) node[midway, opacity=1](C-blue){$C$};
\draw[thick, fill] (0.5,0) circle (0.07);
\draw[thick, fill] (2.5,1) circle (0.07);
\draw[thick, fill] (0,2) circle (0.07);
\draw[thick, fill] (1.5,3) circle (0.07);
\draw[thick, fill] (2,1.5) circle (0.07);
\draw[thick, fill] (3,2.5) circle (0.07);
\draw[thick] (0.5,1) circle (0.07);
\draw[thick] (2.5,0) circle (0.07);
\draw[thick] (0,3) circle (0.07);
\draw[thick] (1.5,2) circle (0.07);
\draw[thick] (2,2.5) circle (0.07);
\draw[thick] (3,1.5) circle (0.07);
\end{scope}

\begin{scope}[xshift=4.5cm]
\draw (-2.5,0) node {$\mathcal{M} \,\, =$};
\draw[thick, fill=black, fill opacity = 0.2]
(-1.4,0) node[left=3pt] (ABC) {} 
-- (-.5, 1.3) node[above=20pt] (BAC) {} node[midway, above=11pt, left=3pt] (Cd) {}
-- (1, 1.3) node[above=20pt] (BCA) {} node[midway, above=20pt] (Bu) {}
-- (1.9, 0) node[right=3pt] (CBA) {} node[midway, above=11pt, right=3pt] (Ad) {}
-- (1,-1.3) node[below=20pt] (CAB) {} node[midway, below=11pt, right=3pt] (Cu) {}
-- (-0.5, -1.3) node[below=20pt] (ACB) {} node[midway, below=20pt] (Bd) {}
-- cycle node[midway, below=11pt, left=3pt] (Au) {};

\begin{scope}[color=blue]
\draw[fill] (-1.4,0) circle (0.05) ;
\draw[fill] (-.5, 1.3) circle (0.05) ;
\draw[fill] (1, 1.3) circle (0.05) ;
\draw[fill] (1.9, 0) circle (0.05) ;
\draw[fill] (1,-1.3) circle (0.05) ;
\draw[fill] (-0.5, -1.3) circle (0.05) ;

\SARdoublybroken{A}{B}{C}
\SARdoublybroken{B}{A}{C}
\SARdoublybroken{B}{C}{A}
\SARdoublybroken{C}{B}{A}
\SARdoublybroken{C}{A}{B}
\SARdoublybroken{A}{C}{B}
\end{scope}

\SARbrokendown{A}{B}{C}
\SARbrokendown{B}{C}{A}
\SARbrokendown{A}{C}{B}

\SARbrokenup{B}{A}{C}
\SARbrokenup{C}{A}{B}
\SARbrokenup{A}{B}{C}

\end{scope}

\end{tikzpicture}
\caption{A domain with Maslov index 3 (on the left), and its moduli space (on the right). The boundary of the moduli space consists of all possible (partially or fully) broken flowlines: in particular, the vertices of the hexagon correspond to doubly broken flowlines, whereas the edges consist of once-broken flowlines.}
\label{fig:ot}
\end{figure}

The same scheme applies if the domain $D$ has Maslov
index $\mu (D)=3$, and $D$ moves 6 coordinates,
so it can be decomposed into three `independent' rectangles, each moving 
different pairs of coordinates, as shown in Figure~\ref{fig:ot}.
The by now customary argument suggests that we need to 
associate to $D$ a permutohedron, i.e. a hexagon in this
dimension. 

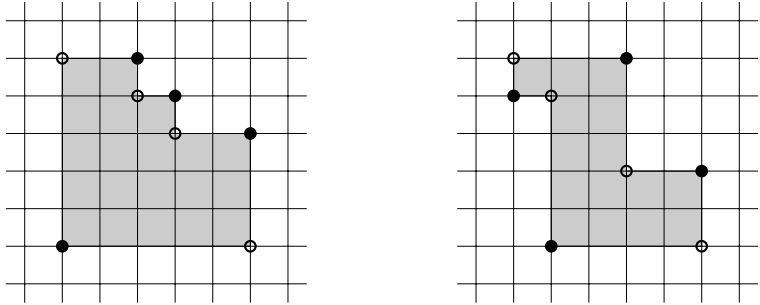
\begin{figure}
\centering
\begin{tikzpicture}

\begin{scope}[xshift=-3cm]
\draw[step=0.5cm,black,very thin] (-0.25,-0.25) grid (3.75,3.75);
\draw[very thin, fill=black, fill opacity=0.2] (0.5, 3) --
(1.5, 3) --
(1.5, 2.5) --
(2, 2.5) --
(2, 2) --
(3,2) --
(3,0.5) --
(0.5,0.5) --
cycle;
\draw[thick, fill] (0.5,0.5) circle (0.07);
\draw[thick, fill] (1.5,3) circle (0.07);
\draw[thick, fill] (2,2.5) circle (0.07);
\draw[thick, fill] (3,2) circle (0.07);
\draw[thick] (0.5, 3) circle (0.07);
\draw[thick] (1.5,2.5) circle (0.07);
\draw[thick] (2,2) circle (0.07);
\draw[thick] (3,0.5) circle (0.07);
\end{scope}

\begin{scope}[xshift=3cm]
\draw[step=0.5cm,black,very thin] (-0.25,-0.25) grid (3.75,3.75);
\draw[very thin, fill=black, fill opacity=0.2] (0.5, 3) --
(2, 3) --
(2, 1.5) --
(3,1.5) --
(3,0.5) --
(1,0.5) --
(1, 2.5) --
(0.5,2.5) --
cycle;
\draw[thick, fill] (0.5,2.5) circle (0.07);
\draw[thick, fill] (2,3) circle (0.07);
\draw[thick, fill] (3,1.5) circle (0.07);
\draw[thick, fill] (1,0.5) circle (0.07);
\draw[thick] (0.5, 3) circle (0.07);
\draw[thick] (2,1.5) circle (0.07);
\draw[thick] (3,0.5) circle (0.07);
\draw[thick] (1,2.5) circle (0.07);
\end{scope}

\end{tikzpicture}
\caption{Other domains with Maslov index 3.}
\label{fig:hat}
\end{figure}

There are many more possibilities for a domain $D$ to have Maslov index equal to 3, and to still have a hexagon as its moduli space $\Mod(D)$; see Figure~\ref{fig:hat} for some examples.

The situation, however, can be much more complicated, and further configurations lead to other natural choices:
the domain shown by Figure~\ref{fig:het},
for example,
provide 4 possible starts from $\x$ when decomposed into rectangles, and altogether 
there are 8 broken trajectories, hence the moduli space is an 
octagon (see Figure \ref{fig:het2}).
\begin{figure}
\centering
\begin{tikzpicture}

\begin{scope}[xshift=-2.5cm]
\draw[step=1cm,black,very thin] (-0.5,-0.5) grid (3.5,3.5);
\draw[very thin, fill=black, fill opacity=0.2]
(0,0) --
(0, 2) --
(1,2) --
(1, 3) --
(3,3) --
(3,1) --
(2,1) --
(2,0);
\draw[very thin, fill=white] (0.1,0.1) rectangle (1.9,0.9);
\draw[very thin, fill=orange, fill opacity=0.2] (0.1,0.1) rectangle (1.9,0.9);
\draw[very thin, fill=white, fill opacity=0.2] (2.9,2.9) rectangle (1.1,2.1);
\draw[very thin, fill=green, fill opacity=0.2] (2.9,2.9) rectangle (1.1,2.1);
\draw[thick, fill] (0,0) circle (0.1);
\draw[thick, fill] (1,2) circle (0.1);
\draw[thick, fill] (3,3) circle (0.1);
\draw[thick, fill] (2,1) circle (0.1);
\draw[thick] (0,2) circle (0.1);
\draw[thick] (1,3) circle (0.1);
\draw[thick] (3,1) circle (0.1);
\draw[thick] (2,0) circle (0.1);
\end{scope}

\begin{scope}[xshift=2.5cm]
\draw[step=1cm,black,very thin] (-0.5,-0.5) grid (3.5,3.5);
\draw[very thin, fill=black, fill opacity=0.2]
(0,0) --
(0, 2) --
(1,2) --
(1, 3) --
(3,3) --
(3,1) --
(2,1) --
(2,0);
\draw[very thin, fill=white] (0.1,0.1) rectangle (0.9,1.9);
\draw[very thin, fill=blue, fill opacity=0.2] (0.1,0.1) rectangle (0.9,1.9);
\draw[very thin, fill=white, fill opacity=0.2] (2.9,2.9) rectangle (2.1,1.1);
\draw[very thin, fill=purple, fill opacity=0.2] (2.9,2.9) rectangle (2.1,1.1);
\draw[thick, fill] (0,0) circle (0.1);
\draw[thick, fill] (1,2) circle (0.1);
\draw[thick, fill] (3,3) circle (0.1);
\draw[thick, fill] (2,1) circle (0.1);
\draw[thick] (0,2) circle (0.1);
\draw[thick] (1,3) circle (0.1);
\draw[thick] (3,1) circle (0.1);
\draw[thick] (2,0) circle (0.1);
\end{scope}

\end{tikzpicture}
\caption{Both sides show a more involved domain with Maslov index 3.
The 4 rectangles shown within the given domain are the 4 possible starts of the 8 decompositions of the domain into rectangles. Each such decomposition corresponds to a fully broken flowlines. 
}
\label{fig:het}
\end{figure}
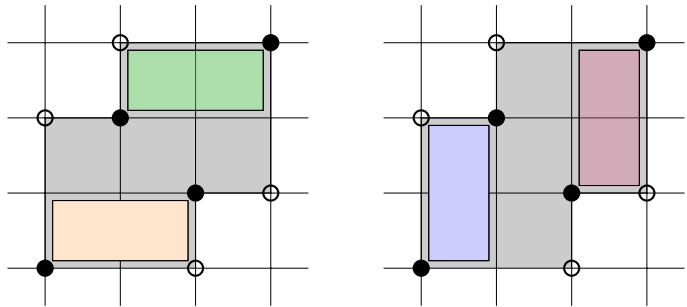
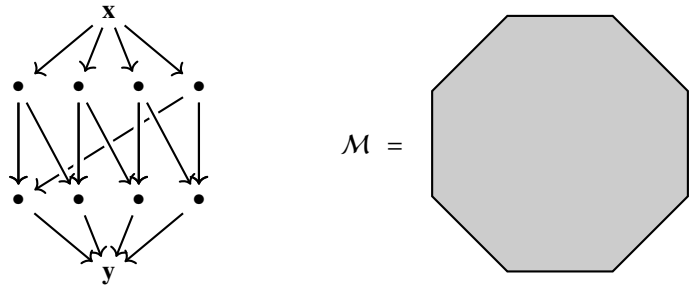
\begin{figure}
\centering
\begin{tikzpicture}[on top/.style={preaction={draw=white,-,line width=#1}},
on top/.default=4pt]

\begin{scope}[xshift=-3cm, yshift=-0.25cm]

\node (x) at (0, 2) {$\x$};

\node (a1) at (-1.2,1) {$\bullet$};
\node (a2) at (-0.4,1) {$\bullet$};
\node (a3) at (0.4,1) {$\bullet$};
\node (a4) at (1.2,1) {$\bullet$};

\node (b1) at (-1.2,-0.5) {$\bullet$};
\node (b2) at (-0.4,-0.5) {$\bullet$};
\node (b3) at (0.4,-0.5) {$\bullet$};
\node (b4) at (1.2,-0.5) {$\bullet$};

\node (y) at (0,-1.5) {$\y$};

\draw[->, thick] (x) -- (a1);
\draw[->, thick] (x) -- (a2);
\draw[->, thick] (x) -- (a3);
\draw[->, thick] (x) -- (a4);

\draw[->, thick] (a4) -- (b1);

\draw[->, thick, on top] (a1) -- (b1);
\draw[->, thick, on top] (a2) -- (b2);
\draw[->, thick, on top] (a3) -- (b3);
\draw[->, thick, on top] (a4) -- (b4);

\draw[->, thick, on top] (a1) -- (b2);
\draw[->, thick, on top] (a2) -- (b3);
\draw[->, thick, on top] (a3) -- (b4);
\draw[->, thick] (a1) -- (b1);
\draw[->, thick] (a2) -- (b2);
\draw[->, thick] (a3) -- (b3);
\draw[->, thick] (a4) -- (b4);

\draw[->, thick] (b1) -- (y);
\draw[->, thick] (b2) -- (y);
\draw[->, thick] (b3) -- (y);
\draw[->, thick] (b4) -- (y);
\end{scope}

\begin{scope}[xshift=3cm]
\def\ml{0.7} 
\def\dp{1} 
\draw (-2.5,0) node {$\mathcal{M} \,\, =$};
\draw[thick, fill=black, fill opacity = 0.2]
(-\ml-\dp,\ml) -- 
(-\ml,\ml+\dp) -- 
(\ml,\ml+\dp) -- 
(\ml+\dp,\ml) -- 
(\ml+\dp,-\ml) -- 
(\ml,-\ml-\dp) -- 
(-\ml,-\ml-\dp) -- 
(-\ml-\dp,-\ml) -- 
cycle;

\end{scope}

\end{tikzpicture}
\caption{The structure of the eight fully broken flowlines and the moduli space associated with the domain of Maslov index 3 from Figure \ref{fig:het}.
The reader is invited to identify the eight fully broken flowlines (corresponding to the eight vertices of the octagon) and the partially broken flowlines (corresponding to the sides of the octagon).}
\label{fig:het2}
\end{figure}

\begin{figure}
\centering
\begin{tikzpicture}[on top/.style={preaction={draw=white,-,line width=#1}},
on top/.default=4pt]

\begin{scope}[xshift=-5cm, yshift=-1cm]
\draw[step=1cm,black,very thin] (-0.5,-0.5) grid (2.5,2.5);
\draw[very thin, fill=black, fill opacity=0.2]
(0,0) rectangle (2,2);
\draw[thick, fill] (0,0) circle (0.1);
\draw[thick, fill] (2,2) circle (0.1);
\draw[thick] (0,2) circle (0.1);
\draw[thick] (2,0) circle (0.1);
\draw[thick, fill] (0.9,1) circle (0.1);
\draw[thick] (1.1,1) circle (0.1);
\end{scope}

\begin{scope}[yshift=-0.25cm]

\node (x) at (0, 2) {$\x$};

\node (a1) at (-0.8,1) {$\bullet$};
\node (a2) at (0.8,1) {$\bullet$};

\node (b1) at (-0.8,-0.5) {$\bullet$};
\node (b2) at (0.8,-0.5) {$\bullet$};

\node (y) at (0,-1.5) {$\y$};

\draw[->, thick] (x) -- (a1);
\draw[->, thick] (x) -- (a2);

\draw[->, thick, on top] (a1) -- (b1);
\draw[->, thick, on top] (a2) -- (b2);

\draw[->, thick, on top] (a1) -- (b2);
\draw[->, thick, on top] (a2) -- (b1);

\draw[->, thick] (b1) -- (y);
\draw[->, thick] (b2) -- (y);
\end{scope}

\begin{scope}[xshift=4.5cm]
\draw (-1.7,0) node {$\mathcal{M} \,\, =$};
\draw[thick, fill=black, fill opacity = 0.2]
(-1,-1) rectangle (1,1);
\end{scope}

\end{tikzpicture}
\caption{A domain with Maslov index 3 (on the left) from $\x$ (full circles) to $\y$ (hollow circles), the fully broken flowlines (in the middle), and the moduli space (on the right). The reader is invited to work out the boundary of the moduli space.}
\label{fig:nyolc}
\end{figure}
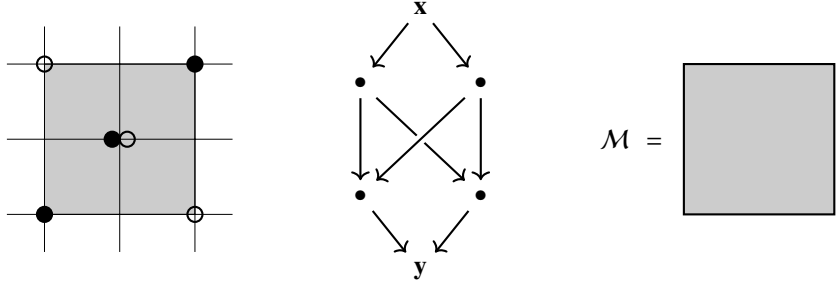

A further example is provided by a rectangle which is not empty
(contains a further coordinate of $\x$ and therefore of $\y$). This
example has Maslov index 3, and there are four broken trajectories in
the decomposition, so the natural choice for the moduli spaces is a
rectangle, depicted in Figure~\ref{fig:nyolc}.

\subsection{Maslov index \texorpdfstring{$k$}{k} domains}
\label{ssec:GeneralCase}
The above \emph{ad hoc} arguments is extended to any
Maslov index in a more systematic way
by induction, using obstruction theory.
There are several complications along the way, though.

First of all, remember that in the flow category we require the moduli spaces to be framed
manifolds (indeed, manifolds with corners). Because of the presence of
`bubbles' in Heegaard Floer homology, however, the natural choices
for spaces associated to (positive)
domains are not always manifolds with corners ---
we need to consider more complicated stratified spaces. On the positive
side, the bubbles will guide us to glue the stratified spaces
associated to the domains connecting two grid states $\x$ and $\y$
to get a framed manifold with corner, which will then be our choice for the
morphism space $\Mod (\x, \y)$ in the flow category.
To explain the main ideas, first we will ignore the complications brought
by the bubbles, and describe the na{\"i}ve inductive idea to define the
manifold with corner associated to a domain $D$ of Maslov index $k+1$.
We will return to the discussion about the role of bubbles in
Subsection~\ref{ssec:bubbles}.

Our inductive hypothesis is therefore that we have already constructed
the framed moduli spaces for all (positive) domains of Maslov index $\leq k$;
consider now a domain $D$ with Maslov index $\mu (D)=k+1$.  By the
inductive hypothesis, the framed boundary $\partial \Mod (D)$ of the
moduli space $\Mod (D)$ is an already defined $(k-1)$-dimensional
framed manifold with corners, and hence provides an element $o(D)$ of the
framed cobordism group $\Omega ^{\mathrm{fr}}_{k-1}$.
The latter is the group of $(k-1)$-dimensional framed
cornered manifolds taken up to framed
cornered cobordism, with group operation given by
the disjoint union.
\begin{remark}
In reality the
boundary defines an element in a more complicated group, which is
denoted by ${\widetilde {\Omega}}^{\mathrm{fr}}_{k-1}$ in
\cite{ManolescuSarkar}: while its definition differs from that of $\Omega ^{\mathrm{fr}}_{k-1}$, it turns out to be isomorphic to it.
Once again, we will ignore this subtlety here.
\end{remark}

Thus, for each (positive) domain $D$ with Maslov index $\mu (D)=k+1$,
we get an element $o(D)\in \Omega^{\mathrm{fr}}_{k-1}$
(the notation indicates that it is an obstruction class).
Putting these classes together, we get a $k$-cochain on a
chain complex generated by the domains, with values in $\Omega
^{\mathrm{fr}}_{k-1}$.

Before proceeding any further, let us discuss the
chain complex generated by the positive domains --- this chain complex
will be called the \emph{obstruction chain complex} and denoted by
$CD_*$.

\begin{definition}
    The chain complex $CD_*=CD_*(\grid )$ is freely generated
    over $\Z$ by the positive domains and graded by their
    Maslov indices:
\[
CD_k(\grid )=\Z \langle (\x,\y,D)\mid D\in \domain ^+(\x,\y), \mu (D)=k\rangle .
\]
Moreover, $CD_*$ is equipped with the endomorphism $\partial$,
given 
by 
\[
\partial (D)=\sum _{
\substack{(R,E)\,\in\, {\mathcal {R}}(\x,\w)\times \domain ^+(\w,\y) \\ R*E=D}
}s(R) E
\,\,\,\,\,+\,\,\,\,\,
(-1)^k \sum _{
\substack{ (E,R)\,\in\, \domain ^+(\x,\w)\times {\mathcal {R}}(\w,\y)\\ E*R=D} 
}s(R)E.
\]
Here ${\mathcal {R}}(\x,\y)$ denotes the rectangles from $\x$
to $\y$, and $s$ denoted a sign assignment.
\end{definition}
It is not hard to see that 
\begin{lemma}
    The pair $(CD _*, \partial )$ is a chain complex, that is, 
    $\partial ^2=0$.
\end{lemma}

On the other hand, it requires a rather tedious calculation to 
show that the homology of this complex is rather simple:
\begin{lemma}
The homology of $(CD_*, \partial)$
is isomorphic to $\Z$, supported in grading 0.
\end{lemma}

Because of the presence of possible bubbles (a.k.a.\ \emph{boundary
  degenerations}), however, later we are forced to consider a more
complicated chain complex, an enhancement $CDP_*$ of $CD_*$ where also
vectors of partitions and further data
describing the positions of the bubbles (and their relative distances)
are added to the generators.  These extra data
will be essential in locating possible bubbles, hence play a crucial
role in gluing arguments.  (We will omit the details of this extension
here.)

 The `earlier defined' cochains (based on the definition of the
 obstruction $o(D)$) naturally live in a subcomplex of the dual $CDP^*$ of $CDP_*$, and an
 important technical step in the construction is the calculation
 showing that this subcomplex is acyclic.


By associating the obstruction 
class $o(D)$ to a Maslov index $k+1$ positive domain,
we get a $(k+1)$-cochain, which turns out to be a cocycle, that is
\[
  \delta (o)=0,
\]
where $\delta$ is the dual coboundary map. This is proved by showing that
$o(\partial D)=0$ for all positive domains $D$ with $\mu(D)=k+2$.
As the cochain complex is acyclic, $o$ is a coboundary, that is,
there is a map $f$ on index $k$ domains associating
\[
E\mapsto f(E)\in \Omega ^{\mathrm{fr}}_{k-1}
\]
to the domain $E$, satisfying 
\[
o=\delta (f).
\]
Now we `change' the framed manifold $\Mod(E)$ (that was already
constructed) to $\Mod(E)\amalg -f(E)$. With respect to these new moduli
spaces, the obstruction class vanishes. In other words, for each
positive domain $D$ with $\mu(D)=k+1$, its boundary $\partial \Mod(D)$
is framed null-cobordant, and therefore, can be filled as a framed
manifold, which we define to be the moduli space $\Mod(D)$. Then the
induction proceeds.  

In this na{\"i}ve approach we would define the moduli spaces $\Mod (\x, \y)$
in the flow category by taking the disjoint union for all
domains of the given Maslov index connecting the two
grid states $\x$ and $\y$.

The potential (and actual) presence of bubbles, however, prevent us to
proceed as outlined above. We need to allow the spaces $\Mod (D)$ to
be more general than framed manifolds with corners; and these
stratified spaces will be glued together in an intricate way to form
$\Mod (\x , \y)$, which will ultimately be a framed manifold with
corners. We will discuss some of the complications caused by the bubbles in
the next subsection through some examples; for the full theory the
interested reader is advised to turn to \cite{ManolescuSarkar}.

\subsection{Bubbles}
\label{ssec:bubbles}

\begin{figure}
\begin{center}
\begingroup%
  \makeatletter%
  \providecommand\color[2][]{%
    \errmessage{(Inkscape) Color is used for the text in Inkscape, but the package 'color.sty' is not loaded}%
    \renewcommand\color[2][]{}%
  }%
  \providecommand\transparent[1]{%
    \errmessage{(Inkscape) Transparency is used (non-zero) for the text in Inkscape, but the package 'transparent.sty' is not loaded}%
    \renewcommand\transparent[1]{}%
  }%
  \providecommand\rotatebox[2]{#2}%
  \newcommand*\fsize{\dimexpr\f@size pt\relax}%
  \newcommand*\lineheight[1]{\fontsize{\fsize}{#1\fsize}\selectfont}%
  \ifx\svgwidth\undefined%
    \setlength{\unitlength}{141.2574708bp}%
    \ifx\svgscale\undefined%
      \relax%
    \else%
      \setlength{\unitlength}{\unitlength * \real{\svgscale}}%
    \fi%
  \else%
    \setlength{\unitlength}{\svgwidth}%
  \fi%
  \global\let\svgwidth\undefined%
  \global\let\svgscale\undefined%
  \makeatother%
  \begin{picture}(1,0.42146997)%
    \lineheight{1}%
    \setlength\tabcolsep{0pt}%
    \put(0,0){\includegraphics[width=\unitlength,page=1]{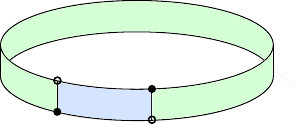}}%
    \put(0.31212494,0.04888889){\color[rgb]{0,0,0}\makebox(0,0)[lt]{\lineheight{1.25}\smash{\begin{tabular}[t]{l}$A$\end{tabular}}}}%
    \put(0.44140854,0.34697355){\color[rgb]{0,0,0}\makebox(0,0)[lt]{\lineheight{1.25}\smash{\begin{tabular}[t]{l}$B$\end{tabular}}}}%
  \end{picture}%
\endgroup%

\end{center}
\caption{Two rectangles forming a width-one strip around the torus.}
\label{fig:kilenc}
\end{figure}

Recall that for two grid states $\x$ and
$\y$ which differ in exactly two coordinates, we associated four
rectangles, two from $\x$ to $\y$ and two from $\y$ to $\x$. Let $A$
be a rectangle from $\x$ to $\y$ and pick a rectangle $B$ from $\y$ to
$\x$.  The union of $A$ and $B$ will comprise a positive domain $D$,
which geometrically is a strip running around the torus (the direction
depending on which rectangle from $\y$ to $\x$ is chosen). See Figure~\ref{fig:kilenc}. If $\mu(D)=2$, then the
strips will have height or width equal to 1 (otherwise they will
contain entire circles, which carry coordinates, driving 
the Maslov index up).

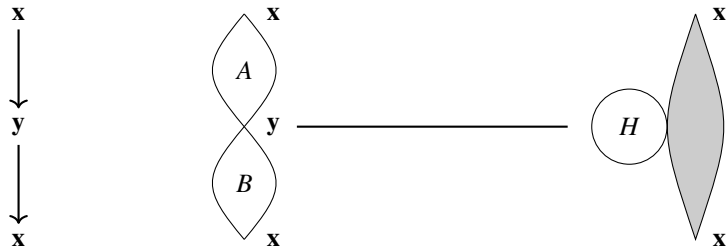
\begin{figure}[h]
\centering
\begin{tikzpicture}[on top/.style={preaction={draw=white,-,line width=#1}},
on top/.default=4pt]

\begin{scope}[xshift=-3cm]

\node (x1) at (0, 1.5) {$\x$};
\node (y) at (0,0) {$\y$};
\node (x2) at (0,-1.5) {$\x$};

\draw[->, thick] (x1) -- (y);
\draw[->, thick] (y) -- (x2);
\end{scope}

\begin{scope}
\draw (0,1.5) node (u1){}; \draw (u1) node[right = 5pt] {$\x$};
\draw (0,0) node (m1){}; \draw (m1) node[right = 5pt] {$\y$};
\draw (0,-1.5) node (t1){}; \draw (t1) node[right = 5pt] {$\x$};
\draw (0,0.75) node {$A$};
\draw (0,-0.75) node {$B$};
\draw (u1.center) edge[out=-130, in=130, looseness=1.5] (m1.center);
\draw (u1.center) edge[out=-50, in=50, looseness=1.5] (m1.center);
\draw (m1.center) edge[out=-130, in=130, looseness=1.5] (t1.center);
\draw (m1.center) edge[out=-50, in=50, looseness=1.5] (t1.center);
\end{scope}

\draw[thick] (0.7,0) -- (4.3,0);

\begin{scope}[xshift=6cm]
\draw (0,1.5) node (u2){}; \draw (u2) node[right = 3pt] {$\x$};
\draw (0,-1.5) node (t2){}; \draw (t2) node[right = 3pt] {$\x$};
\draw[fill=black, fill opacity=0.2] (u2.center) .. controls (-0.5,0) and (-0.5,0) ..
      (t2.center) .. controls (0.5,0) and (0.5,0) .. (u2.center);
\draw (-0.88,0) node {$H$} circle (0.5);
\end{scope}

\end{tikzpicture}
\caption{The other end of this moduli space can be viewed as a boundary degeneration. We follow the convention that horizontal bubbles
will be drawn on the left, while vertical bubbles on the right.}
\label{fig:tiz}
\end{figure}

The composition $D$ of these domains, seen as a domain from $\x$ to $\x$, has a moduli space which is an interval, with one end being the broken flowline given by $A$ (from $\x$ to $\y$) followed by $B$ (from $\y$ to $\x$), and the other end comprising a bubble, see Figure~\ref{fig:tiz}.

\begin{remark}
  In the holomorphic picture the other end of the moduli space, the
  `bubble', corresponds to a geometric picture: it is given by a disk
  which has its entire boundary only on one of the sets of curves (in
  the grid language either only on the horizontal or on the vertical
  ones). This explains the pictorial presentation of
  Figure~\ref{fig:tiz} on the right.
\end{remark}

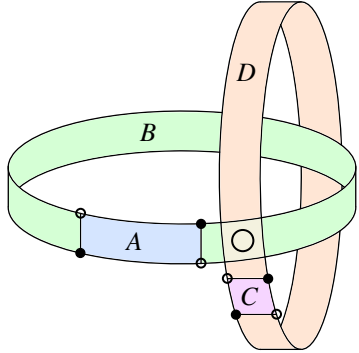
\begin{figure}
\begin{center}
\begingroup%
  \makeatletter%
  \providecommand\color[2][]{%
    \errmessage{(Inkscape) Color is used for the text in Inkscape, but the package 'color.sty' is not loaded}%
    \renewcommand\color[2][]{}%
  }%
  \providecommand\transparent[1]{%
    \errmessage{(Inkscape) Transparency is used (non-zero) for the text in Inkscape, but the package 'transparent.sty' is not loaded}%
    \renewcommand\transparent[1]{}%
  }%
  \providecommand\rotatebox[2]{#2}%
  \newcommand*\fsize{\dimexpr\f@size pt\relax}%
  \newcommand*\lineheight[1]{\fontsize{\fsize}{#1\fsize}\selectfont}%
  \ifx\svgwidth\undefined%
    \setlength{\unitlength}{131.48336936bp}%
    \ifx\svgscale\undefined%
      \relax%
    \else%
      \setlength{\unitlength}{\unitlength * \real{\svgscale}}%
    \fi%
  \else%
    \setlength{\unitlength}{\svgwidth}%
  \fi%
  \global\let\svgwidth\undefined%
  \global\let\svgscale\undefined%
  \makeatother%
  \begin{picture}(1,1.00047124)%
    \lineheight{1}%
    \setlength\tabcolsep{0pt}%
    \put(0,0){\includegraphics[width=\unitlength,page=1]{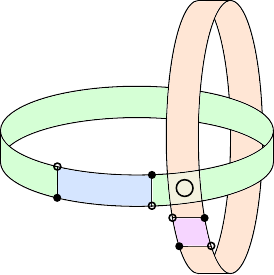}}%
    \put(0.37751907,0.59861056){\color[rgb]{0,0,0}\makebox(0,0)[lt]{\lineheight{1.25}\smash{\begin{tabular}[t]{l}$B$\end{tabular}}}}%
    \put(0.65402714,0.77084645){\color[rgb]{0,0,0}\makebox(0,0)[lt]{\lineheight{1.25}\smash{\begin{tabular}[t]{l}$D$\end{tabular}}}}%
    \put(0.33515917,0.28624233){\color[rgb]{0,0,0}\makebox(0,0)[lt]{\lineheight{1.25}\smash{\begin{tabular}[t]{l}$A$\end{tabular}}}}%
    \put(0.66756606,0.12985113){\color[rgb]{0,0,0}\makebox(0,0)[lt]{\lineheight{1.25}\smash{\begin{tabular}[t]{l}$C$\end{tabular}}}}%
  \end{picture}%
\endgroup%

\end{center}
\caption{The two thin strips around the torus, passing through the fixed $O$-marking. Each can be decomposed into two rectangles, as shown. (The $\x$ coordinates are shown as solid dots.)}
\label{fig:tizenegy}
\end{figure}

The contribution from the bubbles does not create any confusion in grid
\emph{homology}, as these bubbles come in pairs: a horizontal 
bubble will be cancelled by a vertical one.
To exhibit an explicit pairing, we use the $O$-markings: we take the strips passing through a fixed $O$-marking to cancel each other, 
see Figure~\ref{fig:tizenegy}.

\begin{figure}[h]
\centering
\begin{tikzpicture}[on top/.style={preaction={draw=white,-,line width=#1}},
on top/.default=4pt]

\begin{scope}[xshift=-5cm]
\draw (0,1.5) node (u1){}; \draw (u1) node[right = 5pt] {$\x$};
\draw (0,0) node (m1){}; \draw (m1) node[right = 5pt] {$\y$};
\draw (0,-1.5) node (t1){}; \draw (t1) node[right = 5pt] {$\x$};
\draw (0,0.75) node {$A$};
\draw (0,-0.75) node {$B$};
\draw (u1.center) edge[out=-130, in=130, looseness=1.5] (m1.center);
\draw (u1.center) edge[out=-50, in=50, looseness=1.5] (m1.center);
\draw (m1.center) edge[out=-130, in=130, looseness=1.5] (t1.center);
\draw (m1.center) edge[out=-50, in=50, looseness=1.5] (t1.center);

\draw[thick] (0.7,0) -- (2.3,0);

\begin{scope}[xshift=4cm]
\draw (0,1.5) node (u2){}; \draw (u2) node[right = 3pt] {$\x$};
\draw (0,-1.5) node (t2){}; \draw (t2) node[right = 3pt] {$\x$};
\draw[fill=black, fill opacity=0.2] (u2.center) .. controls (-0.5,0) and (-0.5,0) ..
      (t2.center) .. controls (0.5,0) and (0.5,0) .. (u2.center);
\draw (-0.88,0) node {$H$} circle (0.5);
\end{scope}
\end{scope}

\draw (0,0.3) node {glue};
\draw[thick] (-0.3,0) -- (0.3,0);

\begin{scope}[xshift=5cm, xscale=-1]
\draw (0,1.5) node (u1){}; \draw (u1) node[right = 5pt] {$\x$};
\draw (0,0) node (m1){}; \draw (m1) node[right = 5pt] {$\y$};
\draw (0,-1.5) node (t1){}; \draw (t1) node[right = 5pt] {$\x$};
\draw (0,0.75) node {$C$};
\draw (0,-0.75) node {$D$};
\draw (u1.center) edge[out=-130, in=130, looseness=1.5] (m1.center);
\draw (u1.center) edge[out=-50, in=50, looseness=1.5] (m1.center);
\draw (m1.center) edge[out=-130, in=130, looseness=1.5] (t1.center);
\draw (m1.center) edge[out=-50, in=50, looseness=1.5] (t1.center);

\draw[thick] (0.7,0) -- (2.3,0);

\begin{scope}[xshift=4cm]
\draw (0,1.5) node (u2){}; \draw (u2) node[right = 3pt] {$\x$};
\draw (0,-1.5) node (t2){}; \draw (t2) node[right = 3pt] {$\x$};
\draw[fill=black, fill opacity=0.2] (u2.center) .. controls (-0.5,0) and (-0.5,0) ..
      (t2.center) .. controls (0.5,0) and (0.5,0) .. (u2.center);
\draw (-0.88,0) node {$V$} circle (0.5);
\end{scope}
\end{scope}

\end{tikzpicture}
\caption{Gluing the two moduli spaces along the two (vertical and horizontal) boundary degenerations.}
\label{fig:tizenketto}
\end{figure}
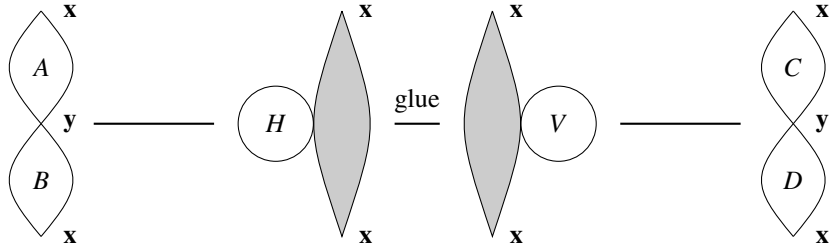

In constructing the flow category, we use these pairs to glue the
moduli spaces $\Mod(A+B)$ and $\Mod(C+D)$ together, along the portion
of their `bubble-type' boundary, as depicted in
Figure~\ref{fig:tizenketto}. In this simple example therefore both
moduli spaces $\Mod (A+B)$ and $\Mod (C+D)$ are intervals, and both
admit two different types of boundaries: one corresponding to a broken
flowline, and the other to a bubble. The gluing is performed on the
bubble-type boundary (where some extra data should be recorded to get
a framed manifold with boundary at the end).
A crucial feature of the end-result is that
the glued-up space no longer has any strata consisting only of bubbles.

The situation encountered above is the first instance when 
bubbles should be taken into account in our constructions.
In grid homology we only need Maslov index 1 domains (to define the boundary
map in the chain complex) and Maslov index 2 domains (to prove that the
square of the boundary map is zero, i.e. we do have a chain complex).
In the spacification process, however, domains with higher Maslov
indices should be also taken into account --- indeed, we need to
work with all positive domains connecting grid states. Therefore
higher dimensional bubbles
should be examined and handled. 

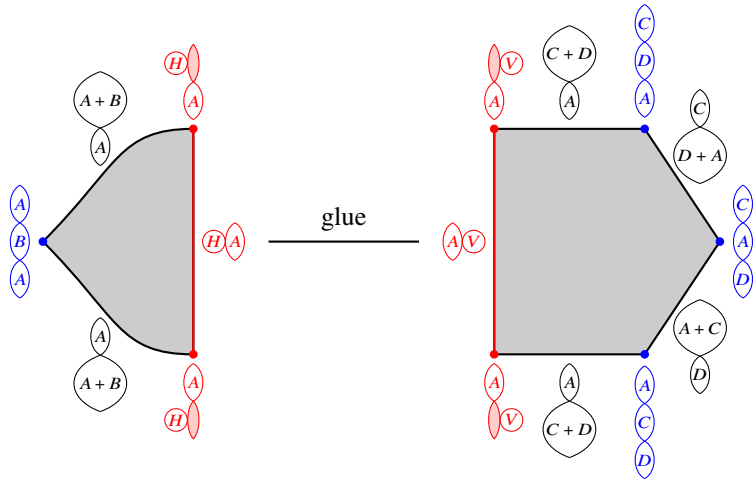
\begin{figure}
\centering
\begin{tikzpicture}

\begin{scope}[xshift=-2cm]
\draw[thick, fill=black, fill opacity = 0.2] (-2,0) node[left=5pt] (ABA){}
.. controls (-1,1) and (-1,1.5) .. node[midway, above=11pt, left=3pt] (Ad){}
(0,1.5) node[right, above=14pt] (HA){}
.. controls (0,0) and (0,0) .. node[midway, right=12pt] (H){}
(0,-1.5) node[right, below=14pt] (AH){}
.. controls (-1,-1.5) and (-1,-1) .. node[midway, below=11pt, left=3pt] (Au){}
(-2,0) -- cycle;
\draw[blue,fill=blue] (-2,0) circle (0.05) ;
\draw[red,fill=red] (0,1.5) circle (0.05) ;
\draw[red,fill=red] (0,-1.5) circle (0.05) ;
\draw[thick, red] (0,-1.5) -- (0,1.5);
\begin{scope}[blue]
\SARdoublybroken{A}{B}{A}
\end{scope}
\SARbrokendown{A}{B}{A}
\SARbrokenup{A}{A}{B}
\begin{scope}[red]
\draw (HA) + (0,0.5) node (HA-T){};
\draw (HA) + (0,-0.5) node (HA-B){};
\draw (HA-T.center) + (-0.1, -0.1) node (HA-T-Bl){};
\draw (HA-T.center) + (0.1, -0.1) node (HA-T-Br){};
\draw (HA.center) + (-0.1, 0.1) node (HA-Bl){};
\draw (HA.center) + (0.1, 0.1) node (HA-Br){};
\draw (HA.center) edge[out=-150, in=150] (HA-B.center);
\draw (HA.center) edge[out=-30, in=30] (HA-B.center);
\draw (HA.center) + (0,-0.25) node {\tiny $A$};
\draw[fill, fill opacity=0.2] (HA-T.center)
.. controls (HA-T-Bl) and (HA-Bl) ..
(HA.center)
.. controls (HA-Br) and (HA-T-Br) ..
(HA-T.center) -- cycle;
\draw (HA.center) + (-0.23,0.25) node {\tiny $H$} circle (0.15);
\draw (AH) + (0,0.5) node (AH-T){};
\draw (AH) + (0,-0.5) node (AH-B){};
\draw (AH.center) + (-0.1, -0.1) node (AH-Bl){};
\draw (AH.center) + (0.1, -0.1) node (AH-Br){};
\draw (AH-B.center) + (-0.1, 0.1) node (AH-B-Bl){};
\draw (AH-B.center) + (0.1, 0.1) node (AH-B-Br){};
\draw (AH-T.center) edge[out=-150, in=150] (AH.center);
\draw (AH-T.center) edge[out=-30, in=30] (AH.center);
\draw (AH.center) + (0,0.25) node {\tiny $A$};
\draw[fill, fill opacity=0.2] (AH.center)
.. controls (AH-Bl) and (AH-B-Bl) ..
(AH-B.center)
.. controls (AH-B-Br) and (AH-Br) ..
(AH.center) -- cycle;
\draw (AH.center) + (-0.23,-0.25) node {\tiny $H$} circle (0.15);
\draw (H) + (0,0.25) node (H-T){};
\draw (H) + (0,-0.25) node (H-B){};
\draw (H-T.center) edge[out=-150, in=150] (H-B.center);
\draw (H-T.center) edge[out=-30, in=30] (H-B.center);
\draw (H.center) node {\tiny $A$};
\draw (H.center) + (-0.28,0) node {\tiny $H$} circle (0.15);
\end{scope}
\end{scope}

\draw (0,0.3) node {glue};
\draw[thick] (-1,0) -- (1,0);

\begin{scope}[xshift=2cm]
\draw[thick, fill=black, fill opacity = 0.2]
(0,1.5) node[left, above=14pt] (VA){}
-- node[midway, above=14pt] (Ad){}
(2, 1.5) node[above=22pt] (CDA){}
-- node[midway, above=22pt, right=3pt] (Cu){}
(3, 0) node[right=5pt] (CAD){}
-- node[midway, below=22pt, right=3pt] (Dd){}
(2, -1.5) node[below=22pt] (ACD){}
-- node[midway, below=14pt] (Au){}
(0,-1.5) node[left, below=14pt] (AV){}
-- node[midway, left=12pt] (V){}
cycle;
\draw[blue,fill=blue] (2,1.5) circle (0.05) ;
\draw[blue,fill=blue] (3,0) circle (0.05) ;
\draw[blue,fill=blue] (2,-1.5) circle (0.05) ;
\draw[red,fill=red] (0,1.5) circle (0.05) ;
\draw[red,fill=red] (0,-1.5) circle (0.05) ;
\draw[thick, red] (0,-1.5) -- (0,1.5);
\begin{scope}[blue]
\SARdoublybroken{C}{D}{A}
\SARdoublybroken{C}{A}{D}
\SARdoublybroken{A}{C}{D}
\end{scope}
\SARbrokendown{C}{D}{A}
\SARbrokenup{C}{D}{A}
\SARbrokendown{A}{C}{D}
\SARbrokenup{A}{C}{D}
\begin{scope}[red]
\draw (VA) + (0,0.5) node (VA-T){};
\draw (VA) + (0,-0.5) node (VA-B){};
\draw (VA-T.center) + (-0.1, -0.1) node (VA-T-Bl){};
\draw (VA-T.center) + (0.1, -0.1) node (VA-T-Br){};
\draw (VA.center) + (-0.1, 0.1) node (VA-Bl){};
\draw (VA.center) + (0.1, 0.1) node (VA-Br){};
\draw (VA.center) edge[out=-150, in=150] (VA-B.center);
\draw (VA.center) edge[out=-30, in=30] (VA-B.center);
\draw (VA.center) + (0,-0.25) node {\tiny $A$};
\draw[fill, fill opacity=0.2] (VA-T.center)
.. controls (VA-T-Bl) and (VA-Bl) ..
(VA.center)
.. controls (VA-Br) and (VA-T-Br) ..
(VA-T.center) -- cycle;
\draw (VA.center) + (0.23,0.25) node {\tiny $V$} circle (0.15);
\draw (AV) + (0,0.5) node (AV-T){};
\draw (AV) + (0,-0.5) node (AV-B){};
\draw (AV.center) + (-0.1, -0.1) node (AV-Bl){};
\draw (AV.center) + (0.1, -0.1) node (AV-Br){};
\draw (AV-B.center) + (-0.1, 0.1) node (AV-B-Bl){};
\draw (AV-B.center) + (0.1, 0.1) node (AV-B-Br){};
\draw (AV-T.center) edge[out=-150, in=150] (AV.center);
\draw (AV-T.center) edge[out=-30, in=30] (AV.center);
\draw (AV.center) + (0,0.25) node {\tiny $A$};
\draw[fill, fill opacity=0.2] (AV.center)
.. controls (AV-Bl) and (AV-B-Bl) ..
(AV-B.center)
.. controls (AV-B-Br) and (AV-Br) ..
(AV.center) -- cycle;
\draw (AV.center) + (0.23,-0.25) node {\tiny $V$} circle (0.15);
\draw (V) + (0,0.25) node (V-T){};
\draw (V) + (0,-0.25) node (V-B){};
\draw (V-T.center) edge[out=-150, in=150] (V-B.center);
\draw (V-T.center) edge[out=-30, in=30] (V-B.center);
\draw (V.center) node {\tiny $A$};
\draw (V.center) + (0.28,0) node {\tiny $V$} circle (0.15);
\end{scope}
\end{scope}

\end{tikzpicture}
\caption{The moduli space $\Mod(2A+B)$ (the triangle) is glued to
    the moduli space $\Mod(A+C+D)$ (the pentagon) along their
    bubble-type boundaries. (Note that the glued-up moduli space is a
    rectangle, not a hexagon.)}
    \label{fig:triangle-pentagon}
\end{figure}

Next we describe a Maslov index $3$
example. Continuing from the previous example, consider the domains $2A+B$ and $A+C+D$ from $\x$ to $\y$. Both the moduli spaces
$\Mod(2A+B)$ and $\Mod(A+C+D)$ should have bubble type boundaries, but
we will glue them carefully so that the glued-up space has no strata consisting only of bubbles.
For $\Mod(2A+B)$, the only
completely broken flowline is $ABA$, and after a little analysis, we
see that $\Mod(2A+B)$ should be a triangle. A similar analysis shows
that $\Mod(A+C+D)$ should be a pentagon. Each of these spaces have one
edge on the boundary corresponding to bubbling, and the edge is
parametrized by the \emph{height} of the bubble.
(We will not discuss the concept of the height of a bubble, or the gluing
parameters appearing later --- these technical data play crucial role
in the gluing theory.)
We then glue the two edges together --- by matching their heights --- and the new moduli space is a smooth manifold with corners,
with no stratum consisting only of bubbles.
This process is illustrated in Figure~\ref{fig:triangle-pentagon}.

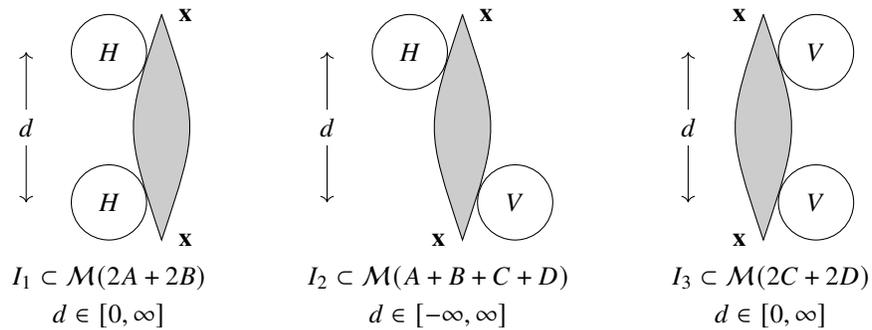
\begin{figure}[h]
\centering
\begin{tikzpicture}

\begin{scope}[xshift=-4cm]
\draw (0,1.5) node (u2){}; \draw (u2) node[right = 3pt] {$\x$};
\draw (0,-1.5) node (t2){}; \draw (t2) node[right = 3pt] {$\x$};
\draw[fill=black, fill opacity=0.2] (u2.center) .. controls (-0.5,0) and (-0.5,0) ..
      (t2.center) .. controls (0.5,0) and (0.5,0) .. (u2.center);
\draw (-0.7,1) node {$H$} circle (0.5);
\draw (-0.7,-1) node {$H$} circle (0.5);
\draw[<->] (-1.8,1) -- node[midway, fill=white] {$d$} (-1.8,-1);
\draw (-0.7,-2) node {$I_1 \subset \mathcal{M}(2A+2B)$};
\draw (-0.7,-2.5) node {$d \in [0, \infty]$};
\end{scope}

\begin{scope}
\draw (0,1.5) node (u2){}; \draw (u2) node[right = 3pt] {$\x$};
\draw (0,-1.5) node (t2){}; \draw (t2) node[left = 3pt] {$\x$};
\draw[fill=black, fill opacity=0.2] (u2.center) .. controls (-0.5,0) and (-0.5,0) ..
      (t2.center) .. controls (0.5,0) and (0.5,0) .. (u2.center);
\draw (-0.7,1) node {$H$} circle (0.5);
\draw (0.7,-1) node {$V$} circle (0.5);
\draw[<->] (-1.8,1) -- node[midway, fill=white] {$d$} (-1.8,-1);
\draw (-0.32,-2) node {$I_2 \subset \mathcal{M}(A+B+C+D)$};
\draw (-0.32,-2.5) node {$d \in [-\infty, \infty]$};
\end{scope}

\begin{scope}[xshift=4cm]
\draw (0,1.5) node (u2){}; \draw (u2) node[left = 3pt] {$\x$};
\draw (0,-1.5) node (t2){}; \draw (t2) node[left = 3pt] {$\x$};
\draw[fill=black, fill opacity=0.2] (u2.center) .. controls (-0.5,0) and (-0.5,0) ..
      (t2.center) .. controls (0.5,0) and (0.5,0) .. (u2.center);
\draw (0.7,1) node {$V$} circle (0.5);
\draw (0.7,-1) node {$V$} circle (0.5);
\draw[<->] (-1,1) -- node[midway, fill=white] {$d$} (-1,-1);
\draw (0.1,-2) node {$I_3 \subset \mathcal{M}(2C+2D)$};
\draw (0.1,-2.5) node {$d \in [0, \infty]$};
\end{scope}

\end{tikzpicture}

\caption{The three intervals, $I_1,I_2,I_3$, of double-bubbles in
    $\Mod(2A+2B)$, $\Mod(A+B+C+D)$, $\Mod(2C+2D)$, respectively. The
    intervals are parametrized by $d$, the height difference between
    the bubbles, and we have $I_1=I_3=[0,\infty]$ and
    $I_2=[-\infty,\infty]$.}
    \label{fig:double-bubble}
\end{figure}

Pushing one dimension further, and continuing from the above examples,
consider the domains $2A+2B$, $A+B+C+D$, and $2C+2D$, all connecting
$\x$ to $\x$. As it was noted earlier, this time the three moduli
spaces $\Mod(2A+2B)$,
$\Mod(A+B+C+D)$, and $\Mod(2C+2D)$ will not be manifolds with
corners, but will be  Whitney stratified
spaces. Each will have bubble-type boundaries, but they will be glued
together carefully to produce a smooth manifold with corners, 
with no stratum consisting only of bubbles.
We will not describe the
full gluing, but rather the gluing near the most complicated portions
of their boundaries, namely the portion of the boundary with two
bubbles. Each of these three moduli spaces has an interval worth of
such `double' bubbling as shown in Figure~\ref{fig:double-bubble}.

\begin{figure}
\definecolor{SARdarkgreen}{RGB}{0,102,0}
\definecolor{SARlightgreen}{RGB}{130,201,85}
\definecolor{SARnavy}{RGB}{0,0,128}
\definecolor{SARlightblue}{RGB}{43,189,255}
\definecolor{SARorange}{RGB}{255,153,85}
  \centering
  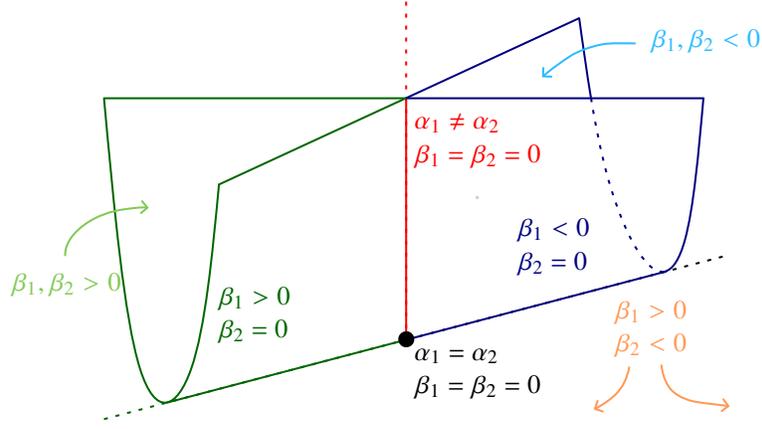
  \caption{The Whitney umbrella stratification of $\R^3$, seen as $\Sym^2(\C)/\R$.
  Let $\alpha_1+i\beta_1$ and
    $\alpha_2+i\beta_2$ denote the two complex numbers, which are
    well-defined up to the $\R$-action
    $(\alpha_1,\alpha_2)\mapsto(\alpha_1+t,\alpha_2+t)$. There is one
    0-dimensional stratum, one 1-dimensional stratum, two
    2-dimensional strata, and three 3-dimensional strata, which are
    shown in different colors and are labelled.}\label{fig:whitney}
\end{figure}

As we got used to gluing horizontal bubbles to
vertical bubbles, we might expect that the `double horizontal
bubbling' of $I_1$ be glued to the `double vertical bubbling' of
$I_3$, and the double `mixed' bubbling of $I_2$ (one horizontal bubble
and one vertical bubble) be glued to itself, hence giving rise to a
quotient of $I_2$ in the gluing.  What really happens is, in fact, more
complicated, because all the above portions of the moduli spaces
corresponding to double bubblings are glued together, as we shall see
in Figure \ref{fig:whitney}. Nonetheless, $I_2$ will be glued to
itself, and this is why in what follows we consider the quotient
$\Mod(A+B+C+D)/\sim$, obtained by quotienting $I_2 \subset \Mod(A+B+C+D)$ by $d\sim -d$.

We now describe how $\Mod(2A+2B)$,
$\Mod(A+B+C+D)/\sim$, and $\Mod(2C+2D)$ are glued in a neighborhood of
$I_1$, $I_2/\sim$, and $I_3$. Each of the intervals $I_1$,
$I_2/\sim$, and $I_3$ is parametrized by $[0,\infty]$, and they are
glued together respecting this parametrization. To extend this gluing
to a neighborhood, we first need to understand the neighborhoods of
these three intervals in $\Mod(2A+2B)$, $\Mod(A+B+C+D)/\sim$, and
$\Mod(2C+2D)$.

Towards that end, consider the stratification of $\R^3$ given by the Whitney umbrella: the ambient space is
$\Sym^2(\C)/\R\cong \R^3$ --- where the $\R$-action is by translating
the real parts --- stratified by the signs of the imaginary
parts. (There is an extra 0-dimensional stratum when the real parts
are equal and both the imaginary parts are 0.) See Figure \ref{fig:whitney}.

Given $(\alpha_1+i\beta_1, \alpha_2+i\beta_2) \in \Sym^2(\C)$, the two real parts
encode the heights of the two bubbles, while the imaginary parts
encode the gluing parameters for gluing the two bubbles, and their
signs encode whether the bubbles are glued in as horizontal bubbles or
as vertical bubbles.

The closure of the red interval (described by the equations $\beta_1=\beta_2=0$) corresponds
to double bubbling, and is parametrized by
$d=|\alpha_1-\alpha_2|\in[0,\infty]$. The interval
$I_2\subset \Mod(A+B+C+D)$ locally looks like
$I_2\subset I_2\times\R_{\geq 0}^2$, where the $\R_{\geq 0}$ are
gluing parameters for gluing the two bubbles. Therefore, after
quotienting, ${I_2/\sim} \subset {\Mod(A+B+C+D)/\sim}$ locally looks like
the red stratum inside the orange stratum. (The orange stratum has one
imaginary part positive and one imaginary part negative, which
corresponds to gluing one bubble in as the horizontal bubble $H$ and
one bubble in as the vertical bubble $V$.)
Similarly,
$I_1\subset \Mod(2A+2B)$ locally looks like the red stratum inside the
light green stratum and $I_3\subset\Mod(2C+2D)$ locally looks like the
red stratum inside the light blue stratum. The Whitney umbrella
witnesses how the light green, orange, and light blue strata glue
together to produce the smooth manifold $\R^3$, and that is the local
model of how the neighborhoods of $I_1,I_2/\sim,I_3$ are glued
together.

\bibliographystyle{alpha}
\bibliography{bibliography}

\end{document}